\RequirePackage{fix-cm}
\documentclass[smallextended]{svjour3}       
\smartqed  
\usepackage{amsmath}
\usepackage{amssymb}
\usepackage{autobreak}
\usepackage[linesnumbered,ruled,vlined]{algorithm2e}
\usepackage{mathrsfs}
\usepackage{graphicx}
\usepackage{subfigure} 
\usepackage{epstopdf}
\usepackage[misc]{ifsym}
\usepackage{bm}
\usepackage{filecontents}
\usepackage{geometry}
\usepackage{makecell}

\usepackage{multicol}
\usepackage{multirow} 
\usepackage{hyperref}
\hypersetup{
	colorlinks=true,
	linkcolor=blue,
	filecolor=gray,      
	urlcolor=blue,
	citecolor=blue,
}

\begin{document}

\title{A robustness-enhanced reconstruction based on discontinuity feedback factor for high-order finite volume scheme 
}


\author{{Hong Zhang}    \textsuperscript{1}   \and
	{Xing Ji}   \textsuperscript{1}
	\and Kun Xu \textsuperscript{2,3,4}
}


\institute{%
	\begin{itemize}
		\item[\textsuperscript{\Letter}] {Xing Ji} \\
		\email{xjiad@connect.ust.hk}
		\at
		\item[\textsuperscript{1}] State Key Laboratory for Strength and Vibration of Mechanical Structures, Shaanxi Key Laboratory of Environment and Control for Flight Vehicle, School of Aerospace Engineering, Xi’an Jiaotong University, Xi’an, China
		\item[\textsuperscript{2}] Department of Mathematics, Hong Kong University of Science and Technology, Clear Water Bay, Kowloon, Hong Kong
		\item[\textsuperscript{3}] Department of Mechanical and Aerospace Engineering, Hong Kong University of Science and Technology, Clear Water Bay, Kowloon, Hong Kong
		\item[\textsuperscript{4}] Shenzhen Research Institute, Hong Kong University of Science and Technology, Shenzhen, China
	\end{itemize}
}

\date{Received: date / Accepted: date}

\authorrunning{Hong Zhang et al.}
\maketitle

\begin{abstract}
	In this paper, a robustness-enhanced reconstruction for the high-order finite volume scheme is constructed on the 2-D structured mesh, and both the high-order gas-kinetic scheme(GKS) and the Lax-Friedrichs(L-F) flux solver are considered to verify the validity of this algorithm. 
	The strategy of the successful WENO reconstruction is adopted to select the smooth sub-stencils. 
	However, there are cases where strong discontinuities exist in all sub-stencils of the WENO reconstruction, which leads to a decrease in the robustness. 
	To improve the robustness of the algorithm in discontinuous regions in two-dimensional space, the hybrid reconstruction based on a combination of discontinuity feedback factor(DF) \cite{ji2021gradient} and WENO reconstruction is developed to deal with the possible discontinuities. 
	Numerical results from smooth to extreme cases have been presented and validate that the new finite volume scheme is effective for robustness enhancement and maintains high resolution compared to the WENO scheme.
\subclass{MSC code1 \and MSC code2 \and more}
\end{abstract}

\section{Introduction}
\label{intro}
In recent decades, great efforts have been devoted to the development of high-order accurate numerical methods for the simulation of compressible flows, with considerable success. The pioneering work in the development of high-order accurate numerical schemes can be attributed to Lax and Wendroff \cite{lax1958systems}. These efforts were further extended to high resolution methods by influential researchers such as Kolgan \cite{kolgan1972application}, van Leer \cite{van1979towards}, Harten \cite{harten1997high} et al. In addition, further advances in high-order numerical methods include variants such as the Essential Non-Oscillatory(ENO) \cite{harten1997uniformly,shu1988efficient}, Weighted Essential Non-Oscillatory(WENO) \cite{liu1994weighted,jiang1996efficient}, and Discontinuous Galerkin(DG) \cite{cockburn1989tvb,cockburn1998runge} etc. 

The GKS is a numerical method for solving the Euler and Navier-Stokes equations, first proposed by Xu \cite{xu1998gas,xu2001gas}. 
The core part of the GKS, i.e., the gas-kinetic flux operates as a Riemann solver, using the Boltzmann-BGK equation to establish the gas distribution function at the cell interface. 
It then calculates the flux by taking into account the relationship between the microscopic gas distribution and the macroscopic physical properties. 
An essential feature of the GKS is the incorporation of a time-dependent gas distribution function at the cell interface, which allows the representation of multi-scale flow phenomena, ranging from kinetic particle transport to hydrodynamic wave propagation. 
The high-order GKS with a single-stage time step method \cite{li2023one} can be achieved by expanding the interface gas distribution function with higher order expressions and spatial components. 
However, beyond the third-order expansion, the use of such high-order GKS is less commonly preferred due to the significantly increased complexity. 
Here we use the GKS based on a two-stage fourth-order(S2O4) temporal discretization \cite{pan2016efficient}. The scheme achieves fourth-order temporal accuracy with only two time advance steps, resulting in a significant reduction in computational cost compared to the fourth-order Runge-Kutta (RK4) method. Subsequently, the S2O4 GKS became the dominant approach in the development of high-order GKS. Subsequent research efforts in the field of high-order GKS have focused mainly on spatial reconstruction.

The spatial reconstruction of high-order GKS is primarily based on the WENO method. 
The WENO method has the advantage of achieving consistent high-order accuracy and high resolution in regions of the solution characterized by smoothness, while simultaneously preserving sharp and essentially monotone shock transitions, and has a wide range of successful applications \cite{borges2008improved,fisher2011boundary,henrick2005mapped}. 
However, the WENO method is too complex and the stencil selection strategy has relatively poor robustness in certain cases, such as high Mach flows of rarefaction waves. 
Recently, there have been many strategies to improve the robustness of the WENO scheme, and the mainstream approach is the hybrid stencils selection strategy, such as by combining two different WENO reconstructions and using a criterion to automatically select the appropriate stencils \cite{liu2013robust}, and the combination of non-linear and linear weights \cite{gao2017enhanced}, which has the advantage of being able to improve the robustness while reducing the amount of computation, but how to choose an appropriate criterion for automatic weight selection is a problem. 
The target ENO provides a new idea to select candidate stencils \cite{fu2016family,fu2018new}. The TENO scheme can automatically reduce the order to a third order depending on the local flow characteristics. This approach effectively handles situations with multiple discontinuities and recovers the robustness of the classical fifth-order WENO scheme \cite{jiang1996efficient,shu1988efficient}. Another effective way to improve the robustness is the limiter, such as the conventional prior limiters, like the van Leer limiter \cite{van1997towards}, the van Albada limiter \cite{van1997comparative}, the Michalak and Ollivier-Gooch limiter \cite{michalak2009accuracy}, the Kitamura–Shima limiter \cite{kitamura2012simple} and posterior limiters like MOOD \cite{clain2011high}, etc. Prior limiters usually have parameters that behave as distinct differences in different cases. Posterior limiters are usually computationally expensive.

In order to improve the robustness of high-order finite volume schemes based on WENO reconstruction, this paper constructs a hybrid spatial reconstruction based on the discontinuity feedback factor (DF) \cite{ji2021gradient} with WENO reconstruction. 
DF is a quantity that uses information about the flow field on either side of a Gaussian point at the cell interface to indicate the discontinuity that may enter the cell at the next time step. 
With DF, the continuous assumption in each cell is no longer necessary for the high-order FVM. 
Therefore, it becomes much more flexible for a high-order FVM to be reduced to a first-order FVM, which can easily achieve the essential positivity preserving property. 
In the present work, DF is used not only for the stencil weights of the hybrid scheme but also for its stencil selection strategy. 
In order to balance the high resolution and robustness of the algorithm, we designed a criterion condition based on the distribution of DF in capturing of strong discontinuities, and the algorithm is approximated to the first-order after satisfying the criterion condition. 
To verify the effectiveness of the algorithm, we tested it on two fluxes. One of the fluxes is the Lax-Friedrichs(L-F) flux function with the strong stability preserving Runge-Kutta(SSP-RK) \cite{gottlieb2005high} temporal discretization, which has the positivity preserving property, and another is the second-order GKS flux function with S2O4 temporal discretization. The robustness of the hybrid reconstruction of the two schemes will be presented.
This paper is organized as follows: Sect.~\ref{sec:1} provides a brief overview of the high-order finite volume scheme. Sect.~\ref{sec:4} introduces the new hybrid reconstruction based on a combination of DF and WENO reconstruction. In Sect.~\ref{sec:5}, we present the results of numerical tests showing the performance and robustness of our improved finite volume scheme in different flow scenarios, and the last section is the conclusion.

\section{High-order finite volume scheme}
\label{sec:1}
\subsection{Finite volume framework}

The Finite Volume Method (FVM) divides the computational domain into a finite number of control volumes, within each control volume, the physical quantities are integrated and averaged, resulting in the discrete form of the fluid dynamics partial differential equations
\begin{equation}
	\int_{{\rm \Omega}_{ij}}\mathbf{W}(\mathbf{x}, t^{n+1}){\rm d}V=\int_{{\rm \Omega}_{ij}}\mathbf{W}(\mathbf{x}, t^n){\rm d}V-\int^{t^{n+1}}_{t_n}\int_{\partial {\rm \Omega}_{ij}}\mathbf{F}\cdot\mathbf{n}{\rm d}S.
\end{equation}
where $\mathbf{W}=(\rho, \rho U, \rho V, \rho E)^T$ are the conservative flow variables in a control volume ${\rm \Omega}_{ij}$, $\partial {\rm \Omega}_{ij}$ corresponds to the cell interface, and $\mathbf{F}$ are the corresponding fluxes across the $\partial{\rm \Omega}_{ij}$. In a 2-D rectangular mesh, the boundary can be expressed as
\begin{equation}
	\partial{\rm \Omega}_{ij}=\bigcup^4_{p=1}{\rm\Gamma}_{ij,p}.
\end{equation}
where ${\rm\Gamma}_{ij,p}$ denotes the $p$th interface of the cell ${ij}$. Integrating over the cell ${\rm \Omega}_{ij}$, the semi-discrete form can be obtained as follows \cite{ji2020hweno}
\begin{equation}\label{eq2}
	\frac{{\rm d}\mathbf{W}_{ij}}{{\rm d}t}=\mathcal{L}(\mathbf{W}_{ij})=-\frac1{\left|{\rm \Omega}_{ij}\right|}\sum^4_{p=1}\oint_{\Gamma_{ij,p}}\mathbf{F}(\mathbf{W}_{ij})\cdot\mathbf{n}_p{\rm d}s.
\end{equation}
where $\mathbf{W}_{ij}$ is the cell averaged value over the cell ${\rm\Omega}_{ij}$, $\left|{\rm \Omega}_{ij}\right|$ is the area of ${\rm \Omega}_{ij}$, $\mathcal{L}(\mathbf{W})$ are the spatial derivatives of the flux, $\mathbf{F}=(F, G)^T$ and $\mathbf{n}_p$ is the outer normal direction of the interface ${\rm \Gamma}_{ij,p}$.

In this paper, we propose a hybrid fifth-order spatial reconstruction method, and both Gaussian quadrature weights $A_1, A_2$ are equal to 1/2 when two Gaussian points are used, and the line integral over ${\rm \Gamma}_{ij,p}$ is discretized according to Gaussian quadrature as follows
\begin{equation}
	\oint_{\Gamma_{ij,p}}\mathbf{F}(\mathbf{W}_{ij})\cdot\mathbf{n}_p{\rm d}s=\left|l_p\right|\sum^2_{m=1}A_m \mathbf{F}(\mathbf{x}_{p,m}, t)\cdot\mathbf{n}_p
\end{equation}
where $\left|l_p\right|$ is the length of the cell interface ${\rm \Gamma}_{ij,p}$, and $\mathbf{x}_{p,m},\ m=1,2$ for ${\rm \Gamma}_{ij,p}$ are the Gaussian points. For the line with endpoints $\mathbf{X}_{p,1}$ and $\mathbf{X}_{p,2}$, the Gaussian quadratures are $\mathbf{x}_{p,1}=c\mathbf{X}_{p,1}+(1-c)\mathbf{X}_{p,2},\mathbf{x}_{p,2}=(1-c)\mathbf{X}_{p,1}+c\mathbf{X}_{p,2}$, where $c=\frac{1}{2}+\frac{\sqrt{3}}{6}$.

To update the flow variables in global coordinates, we need to calculate the fluxes in unit length across each Gaussian point as follows
\begin{equation}
	\mathbf{F}(\mathbf{x}_{p,m}, t)\cdot \mathbf{n}_p=\int \bm{\psi} f(\mathbf{x}_{p,m},t,\mathbf{u},\xi)\mathbf{u}\cdot\mathbf{n}{\rm d}\mathbf{u}{\rm d}\xi
\end{equation}
where $f(\mathbf{x}_{p,m},t,\mathbf{u},\xi)$ is the gas distribution function at the corresponding Gaussian point. The rotation matrix $\mathbf{T}$ from global to local coordinates for the 2-D case is
$$
\mathbf{T}=
\setlength{\arraycolsep}{8pt}
\left(
\begin{array}{cccc}
	1 & 0 & 0 & 0 \\
	0 & {\rm cos}\,\theta & {\rm sin}\,\theta & 0 \\
	0 & -{\rm sin}\,\theta & {\rm cos}\,\theta & 0 \\
	0 & 0 & 0 & 1 \\
\end{array}
\right).
$$
then we can obtain the conservative variables in the local coordinate
\begin{equation}
	\mathbf{\widetilde{W}}=\mathbf{TW}.
\end{equation}
where $\mathbf{\widetilde{W}}=(\tilde{\rho}, \tilde{\rho}\tilde{u}, \tilde{\rho}\tilde{v}, \tilde{\rho}\tilde{E})$. Here
we first calculate the fluxes in the local coordinate
\begin{equation}
	\mathbf{\tilde{F}}_{p,m}(t)=\int\bm{\psi}\tilde{u}\tilde{f}(\tilde{\mathbf{x}}_{p,m},t,\tilde{\mathbf{u}}, \xi){\rm d}\tilde{\mathbf{u}}{\rm d}\xi.
\end{equation}
where the initial point of the local coordinate is $\tilde{\mathbf{x}}_{p,m}=(0,0)$ with the x direction in $\mathbf{n}_p$.
Finally we can obtain the fluxes in the global coordinate. In the 2-D case, the global and local fluxes are related as \cite{pan2016third}
\begin{equation}
	\mathbf{F}(\mathbf{W})\cdot \mathbf{n}=\mathbf{T}^{-1}\widetilde{\mathbf{F}}(\mathbf{\widetilde{W}}).
\end{equation}

\subsection{Gas-kinetic flux solver and two-stage fourth-order temporal discrization with non-linear time limiter}
The 2-D BGK equation \cite{1} can be written as
\begin{equation}\label{eq1}
	f_t+{\mathbf{u}}\cdot \nabla f=\frac{g-f}{\tau}.
\end{equation}
where $f$ is the gas distribution function, $g$ is the corresponding equilibrium state, and $\tau$ is correspond to the collision time. The collision term satisfies the compatibility condition
\begin{equation}
	\int\frac{g-f}{\tau}\bm{\psi} d{\rm \Xi}=0.
\end{equation} 
where $\bm{\psi} =\left(1, u, v, \frac12\left(u^2+v^2+\xi^2\right)\right)^T,\ d{\rm \Xi}=dudvd\xi_1\cdots d\xi_K$, $(u,v)$ are the two components of the macroscopic particle velocity, $\xi_i=(\xi_1,\xi_2,\cdots,\xi_K)$ are the components of the internal particle velocity in $K$ dimensions.
$K$ is the number of internal freedom. $K = (4-2\gamma)/(\gamma-1)$ for 2-D flows, and $\gamma$ is the specific heat ratio.

The general macroscopic gas dynamic equations can be obtained by expanding the BGK equation by the Chapman-Enskog expansion \cite{chapman1990mathematical}, where the gas distribution function can be expressed as
\begin{equation}
	f=g-\tau D_{\mathbf{u}}g+\tau D_{\mathbf{u}}(\tau D_{\mathbf{u}})g-\tau D_{\mathbf{u}}[\tau D_{\mathbf{u}}(\tau D_{\mathbf{u}})g]+\cdots,
\end{equation}
where $D_{\mathbf{u}}=\partial/\partial t+\mathbf{u}\cdot\nabla$. For example, $f=g$ corresponds to the Euler equations and the Navier-Stokes equations can be obtained by the truncated first-order distribution function
\begin{equation}
	f=g-\tau(ug_x+vg_y+g_t).
\end{equation}

The key to the GKS is the construction of the time-dependent cell interface distribution function $f$. The integral solution of the BGK equations(\ref{eq1}) is \cite{xu2014direct}
\begin{equation}
	f(\mathbf{x}_{p,m},t,\mathbf{u},\xi)=\frac{1}{\tau}\int^t_0g(\mathbf{x}^{\prime} ,t^{\prime},\mathbf{u},\xi)e^{-(t-t^{\prime})/\tau}{\rm d}t^{\prime}+e^{-t/\tau}f_0(\mathbf{x}_{p,m}-\mathbf{u}t, \mathbf{u}, \xi).
\end{equation}
where $\mathbf{x}_{p,m}=(0,0)$ is the quadrature point at the interface in the local coordinates, and $\mathbf{x}_{p,m}=\mathbf{x}^{\prime}+\mathbf{u}(t-t^{\prime})$ is the trajectory of the particles, $f_0$ is the initial gas distribution function, and $g$ is the corresponding equilibrium state. Then we can construct a second-order time accurate gas distribution function \cite{xu2001gas} at the local Gaussian point $\mathbf{x}_{p,m}=(0,0)$
\begin{equation}
	\begin{aligned}
		&f(\mathbf{x}_{p,m},t,\mathbf{u},\xi)=\left(1-e^{-t/\tau_n}\right)g^c\\
		&+\left[(t+\tau)e^{-t/\tau_n}-\tau\right]a^c_\mathbf{x}\cdot\mathbf{u}g^c+\left(t-\tau+\tau e^{-t/\tau_n}\right)A^cg^c\\
		&+e^{t/\tau_n}g^l\left[1-(\tau+t)a^l_\mathbf{x}\cdot\mathbf{u}-\tau A^l\right]\mathbb{H}(u)\\
		&+e^{t/\tau_n}g^r\left[1-(\tau+t)a^r_\mathbf{x}\cdot\mathbf{u}-\tau A^r\right]\left[1-\mathbb{H}(u)\right].\\
	\end{aligned}
\end{equation}
where $\mathbb{H}$ is the Heaviside function. The function $g^{l,r,c}$ represents the initial gas distribution function at the left and right sides of the cell interface, and the equilibrium state located at an interface, respectively. The flow behavior at the cell interface depends on the ratio of the time step to the local particle collision time $\Delta t/\tau$.

The function $g^k,\ k=l,r$ satisfies Maxwell's distribution
\begin{equation}
	g^k=\rho^k\left(\frac{\lambda^k}{\pi}\right)e^{-\lambda^k\left[\left(u-U^k\right)^2+\left(v-V^k\right)^2+\xi^2\right]}.
\end{equation}
where $\lambda$ is a function of temperature, molecular mass and the Boltzmann constant. $g^k$ can be determined by the macroscopic variables $\mathbf{W}^l, \mathbf{W}^r$ through spatial reconstruction(see Sect.~\ref{sec:4})
\begin{equation}
	\int \bm{\psi}g^l{\rm d}\Xi=\mathbf{W}^l,\quad \int \bm{\psi}g^r{\rm d}\Xi=\mathbf{W}^r.
\end{equation}
The coefficients $a_{\mathbf{x}},A$ denote the spatial and temporal derivatives, respectively, which have the form
\begin{equation}
	a_{\mathbf{x}}\equiv(\partial g/\partial \mathbf{x})/g=g_{\mathbf{x}}/g,\quad A\equiv (\partial g/\partial t)/g=g_t/g.
\end{equation}
which are determined by the spatial derivatives of $\mathbf{W}$ and the compatibility condition as follows
\begin{equation}
	\begin{aligned}
		&\left<a_x\right>=\frac{\partial \mathbf{W}}{\partial x}=\mathbf{W}_x,\quad \left<a_y\right>=\frac{\partial\mathbf{W}}{\partial y}=\mathbf{W}_{y},\\
		&\left<A+a_xu+a_yv\right>=0.
	\end{aligned}
\end{equation}
where $a_x=(a_{x1},\ a_{x2}u,a_{x3}v,\ a_{x4}\frac12(u^2+v^2+\xi^2))^T$ and $a_y$ has the similar form. $\left<\cdots\right>$ are the moments of a gas distribution function defined by
\begin{equation}
	\left<(\cdots)\right>=\int\bm{\psi}(\cdots)g{\rm d}\Xi.
\end{equation}
Similarly, the equilibrium state $g^c$ and its derivatives $a^c_{\mathbf{x}},A^c_{\mathbf{x}}$ are determined by the corresponding $\mathbf{W}^c,\mathbf{W}^c_{\mathbf{x}}$, and the construction of the $\mathbf{W}^c,\mathbf{W}^c_{\mathbf{x}}$ can be found in Sect.~\ref{sec:4}. More details about the integration calculation can be found in Ref \cite{ji2019high}. 

The physical collision time $\tau$ in the exponential function can be replaced by a numerical collision time $\tau_n$ for capture the unresolved discontinuity. For the inviscid flow, the collision time $\tau_n$ has the form \cite{xu2001gas}
\begin{equation}
	\tau_n=C_1\Delta t +C_2\left|\frac{p_l-p_r}{p_l+p_r}\right|\Delta t.
\end{equation}
where $C_1=0.01, C_2=5.0$ are used in this paper. When $C_2=5.0$, the maximum coefficient on the left side of $\Delta t$ approaches 5, and $e^{-1/5}$ approaches 1.0(0.8187), such values are advantageous for the robustness of the numerical examples in this paper. $p_l,p_r$ denote the pressure on the left and right sides of the cell interface, respectively. For the viscous flow, see Ref \cite{xu2001gas}.

The two-stage fourth-order(S2O4) temporal discretization was developed for L-F solvers, and has been used in CFD applications \cite{li2016two,pan2016efficient}, which has the form
\begin{equation}\label{eq24}
	\begin{aligned}
		&\mathbf{W}^*_{ij}=\mathbf{W}^n_{ij}+\frac12\Delta t\mathcal{L}(\mathbf{W}_{ij}^n)+\frac18\Delta t^2\frac{\partial}{\partial t}\mathcal{L}(\mathbf{W}_{ij}^n),\\
		&\mathbf{W}^{n+1}_{ij}=\mathbf{W}_{ij}^n+\Delta t\mathcal{L}(\mathbf{W}^n_{ij})+\frac{1}{6}\Delta t^2\left(\frac{\partial}{\partial t}\mathcal{L}(\mathbf{W}^n_{ij})+2\frac{\partial}{\partial t}\mathcal{L}(\mathbf{W}^*_{ij})\right).
	\end{aligned}
\end{equation}
where $\frac{\partial}{\partial t}\mathcal{L}(\mathbf{W})$ are the time derivatives of the summation of the flux transport at the closed interface of the cell.

For the GKS, the flux function is a time-dependent complicated function, and to evaluate the flux function, it is necessary to obtain the first-order time derivatives of the flux at $t_n$ and $t_*=t_n+\Delta t/2$. The total flux transport at the Gaussian point $\mathbf{x}_{p,m}$ over a time interval $\delta$ is denoted by
\begin{equation}\label{eq25}
	\mathbb{F}_{p,m}(\mathbf{W}^n,\delta)=\int_{t_n}^{t_n+\delta}\mathbf{F}_{p,m}(\mathbf{W}^n,t)dt=\int_{t_n}^{t_n+\delta}\int u\bm{\psi}f(\mathbf{x}_{p,m},t,\mathbf{u},\xi){\rm d}\Xi{\rm d}t.
\end{equation}

For second-order gas-kinetic solver, assuming $t_n=0$, the flux can be approximated as a linear function
\begin{equation}\label{eq26}
	\mathbf{F}_{p,m}(\mathbf{W}^n,t)=\mathbf{F}_{p,m}^n+\partial_t\mathbf{F}^n_{p,m}t.
\end{equation}
Substituting Eq. (\ref{eq26}) into Eq. (\ref{eq25}), the coefficients $\mathbf{F}_{p,m}^n$ and $\partial_t\mathbf{F}^n_{p,m}$ can be determined as follows
\begin{equation}
	\begin{aligned}
		&\mathbf{F}^n_{p,m}\Delta t+\frac12\partial_t\mathbf{F}_{p,m}^n\Delta t^2=\mathbb{F}(\mathbf{W}^n,\Delta t),\\
		&\frac12\partial_t\mathbf{F}_{p,m}^n\Delta t+\frac18\partial_t\mathbf{F}_{p,m}^n\Delta t^2=\mathbb{F}(\mathbf{W}^n,\Delta t/2).
	\end{aligned}
\end{equation}
By solving the linear system, we have
\begin{equation}
	\begin{aligned}
		&\mathbf{F}_{p,m}^n=(4\mathbb{F}_{p,m}(\mathbf{W}^n, \Delta t/2)-\mathbb{F}_{p,m}(\mathbf{W}^n, \Delta t))/\Delta t,\\
		&\partial_t\mathbf{F}_{p,m}^n=4(\mathbb{F}_{p,m}(\mathbf{W}^n, \Delta t)-2\mathbb{F}_{p,m}(\mathbf{W}^n, \Delta t/2))/\Delta t^2.
	\end{aligned}
\end{equation}
then $\mathcal{L}(\mathbf{W}^n_{ij})$ and $\frac{\partial}{\partial_t}\mathcal{L}(\mathbf{W}^n_{ij})$ can be obtained
\begin{equation}
	\begin{aligned}
		&\mathcal{L}(\mathbf{W}^n_{ij})=-\frac{1}{\left|{\rm \Omega}_{ij}\right|}\sum^{4}_{p=1}\sum^2_{m=1}A_m\mathbf{F}(\mathbf{x}_{p,m},t_n)\cdot \mathbf{n}_p,\\
		&\frac{\partial}{\partial_t}\mathcal{L}(\mathbf{W}^n_{ij})=-\frac{1}{\left|{\rm \Omega}_{ij}\right|}\sum^{4}_{p=1}\sum^2_{m=1}A_m\partial_t\mathbf{F}(\mathbf{x}_{p,m},t_n)\cdot \mathbf{n}_p.\\
	\end{aligned}
\end{equation}
According to Eq. (\ref{eq24}), $\mathbf{W}^*_i$ can be obtained, and a similar formulation for $\mathbf{W}^{n+1}_{ij}$. 

To validate the effectiveness of the robustness-enhancement technique, the GKS solver with S2O4 time discretization is considered. However, the possible discontinuities of the flux function in time have not been considered in the S2O4 method, which will affect the evaluation of the stability of the algorithm. To limit the possible discontinuities in time, the non-linear time discretization limiter is used \cite{zhao2023direct} and the S2O4 method can be reformulated as
\begin{equation}
	\mathbf{W}^{n+1}_{ij}=\mathbf{W}^n_{ij}+\Delta t \mathcal{L}(\mathbf{W}^n_{ij})+\frac{\Delta t^2}2\partial_t\mathcal{L}(\mathbf{W}^n_{ij})-\frac{\Delta t^2}{3}\partial_t\tilde{\mathcal{L}}(\mathbf{W}^n_{ij})+\frac{\Delta t^2}{3}\partial_t\tilde{\mathcal{L}}(\mathbf{W}^*_{ij}).
\end{equation}
where $\partial_t\tilde{\mathcal{L}}(\mathbf{W}^n_{ij})$ and $\partial_t\tilde{\mathcal{L}}(\mathbf{W}^*_{ij})$ are the limited time derivatives. For example, based on $\partial_t{\mathcal{L}}(\mathbf{W}^n_{ij})$, $\partial_t\tilde{\mathcal{L}}(\mathbf{W}^n_{ij})$ has the form
\begin{equation}
	\partial_t\tilde{\mathcal{L}}(\mathbf{W}^n_{ij})=-\frac1{{\rm \Omega}_{ij}}\sum^4_{p=1}\omega_p^t\left(\sum^2_{m=1}A_m\partial_t\mathbf{F}(\mathbf{x}_{p,m},t_n)\right)\cdot \mathbf{n}_p.
\end{equation}
where $\omega^t_p$ is a nonlinear weight for the $p$th interface of the cell, which can be obtained by
\begin{equation}
	\begin{aligned}
		\tilde{\alpha}^k_1&=1+\left(\frac{\tau^k_Z}{\beta^k_{\rm min}+\varepsilon}\right)^2,\ \tilde{\alpha}^k_2=1+\left(\frac{\tau^k_Z}{\beta^k_{\rm max}+\varepsilon}\right)^2,\alpha^k=2\frac{\tilde{\alpha}^k_2}{\tilde{\alpha}^k_1+\tilde{\alpha}^k_2}\ k=L,R\\
		\omega^t_p&={\rm min}\{\alpha^L,\alpha^R\}.
	\end{aligned}
\end{equation}
where $\varepsilon$ is a small value and $\varepsilon=1\times10^{-6}$ is used in this paper. $\beta^{L,R}_{\rm min},\beta^{L,R}_{\rm max},\tau^{L,R}_Z$ are used as smoothness indicators to evaluate the smoothness of the distribution of flow variables near the cell interface. The definition of these smoothness indicators is given in Eq.~(\ref{eq37}) and Eq.~(\ref{eq38}) in the Sect~\ref{sec:3.2}.
More details on non-linear limiter for temporal discretization can be found in Ref \cite{zhao2023direct}.

\subsection{Lax-Friedrichs flux solver and RK temporal discrization}
The 2-D Euler equations can be written as
\begin{equation}
	\left(
	\begin{array}{c}
		\rho \\
		\rho U \\
		\rho V\\
		\rho E
	\end{array}
	\right)_t  +
	\left(
	\begin{array}{c}
		\rho U\\
		\rho U^2+p \\
		\rho UV\\
		U(\rho E+p)
	\end{array}
	\right)_x+
	\left(
	\begin{array}{c}
		\rho V\\
		\rho UV\\
		\rho V^2+p \\
		V(\rho E+p)
	\end{array}
	\right)_y=0
\end{equation}
The Lax-Friedrichs(L-F) method \cite{toro2013riemann} is a numerical method for the solution of hyperbolic partial differential equations based on finite differences. The method can be described as the forward-in-time, centered-in-space (FTCS) scheme with a numerical dissipation term of 1/2. The first-order L-F flux function can be expressed as 
\begin{equation}
	\begin{aligned}
		\mathbf{F}(\mathbf{W}_{i+1/2})&=\frac12(\mathbf{F}(\mathbf{W}^l_{i+1/2})+\mathbf{F}(\mathbf{W}^r_{i+1/2}))\\
		&-\frac12{\rm max}\{|U^l_{i+1/2}|+c^l_{i+1/2},|U^r_{i+1/2}|+c^r_{i+1/2}\}(\mathbf{W}^r_{i+1/2}-\mathbf{W}^l_{i+1/2})
	\end{aligned}
\end{equation}
where $\mathbf{F}(\mathbf{W}^l)=\left(\rho^lU^l, \rho^lU^lU^l+p^l,\rho^lU^lV^l,U^l(\rho^lE^l+p^l)\right)^T$ corresponds to the flux on the left side of the interface, and $c^l=\sqrt{\frac{\gamma p^l}{\rho^l}}$ is the speed of sound.

The third-order strong stability preserving Runge-Kutta(SSP-RK3) \cite{gottlieb2005high} time discretization has been developed for the solution of semi-discrete method of lines approximations of hyperbolic partial differential equations. The SSP-RK3 method preserves the strong stability properties in any norm or seminorm of the spatial discretization coupled with a first-order Euler time step of the form
\begin{equation}\label{eq24}
	\begin{aligned}
		&\mathbf{W}^{(1)}_{ij}=\mathbf{W}^n_{ij}+\Delta t\mathcal{L}(\mathbf{W}_{ij}^n),\\
		&\mathbf{W}^{(2)}_{ij}=\frac34\mathbf{W}^n_{ij}+\frac14\mathbf{W}^{(1)}_{ij}+\frac14\Delta t\mathcal{L}(\mathbf{W}_{ij}^{(1)}),\\
		&\mathbf{W}^{n+1}_{ij}=\frac13\mathbf{W}^n_{ij}+\frac23\mathbf{W}^{(2)}_{ij}+\frac23\Delta t\mathcal{L}(\mathbf{W}_{ij}^{(2)}).\\
	\end{aligned}
\end{equation}

The L-F solver with the SSP-RK method has the positive-preserving property, which can guarantee the robustness of the algorithm, which is only affected by the spatial reconstruction.

\section{New hybrid reconstruction based on discontinuity feedback factor and WENO}
\label{sec:4}
In this section, we have developed a hybrid spatial reconstruction method based on DF and WENO techniques, which significantly increases robustness while maintaining high resolution.
\subsection{WENO-AO reconstruction}\label{sec:3.2}
The key idea of the WENO-AO method is to reconstruct a reliable polynomial based on three sub-stencils and a large-stencil \cite{balsara2016efficient}, which can obtain fifth-order accuracy in the smooth region.

In the 2-D case, the reconstruction is performed direction by direction. For example, at the cell interface$(i+1/2,j)$, the 1-D WENO-Z reconstruction \cite{borges2008improved} is first applied to obtain the interface value $\mathbf{Q}^l_{i+1/2,j}$ by using the cell values $\mathbf{Q}_{i-2},\cdots,\mathbf{Q}_{i+2}$, where $\mathbf{Q}$ are the characteristic variables calculated by $\mathbf{W}$. Then, the tangential reconstruction is performed by using the 1-D WENO-Z again in the y-direction to obtain the values at the Gaussian points based on $\mathbf{Q}^l_{i+1/2,j-2},\cdots\mathbf{Q}^l_{i+1/2,j+2}$. 

For example, the normal reconstruction in the x-direction, $(i+1/2,j)$ corresponds to $(0,0)$, and the reconstruction polynomials form are as follows
\begin{equation}
	\begin{aligned}
		p^{r3}_0(x)&=\frac13\mathbf{Q}_{i-2,j}-\frac76\mathbf{Q}_{i-1,j}+\frac{11}{6}\mathbf{Q}_{i,j}+\frac1{\Delta x}(\mathbf{Q}_{i-2,j}-3\mathbf{Q}_{i-1,j}+2\mathbf{Q}_{i,j})x\\
		&+\frac1{\Delta x^2}\left(\frac12\mathbf{Q}_{i-2,j}-\mathbf{Q}_{i-1,j}+\frac12\mathbf{Q}_{i,j}\right)x^2,\\
		p^{r3}_1(x)&=-\frac16\mathbf{Q}_{i-1,j}+\frac56\mathbf{Q}_{i,j}+\frac{1}{3}\mathbf{Q}_{i+1,j}+\frac1{\Delta x}(\mathbf{Q}_{i+1,j}-\mathbf{Q}_{i,j})x\\
		&+\frac1{\Delta x^2}\left(\frac12\mathbf{Q}_{i-1,j}-\mathbf{Q}_{i,j}+\frac12\mathbf{Q}_{i+1,j}\right)x^2,\\
		p^{r3}_2(x)&=\frac13\mathbf{Q}_{i,j}+\frac56\mathbf{Q}_{i+1,j}-\frac{1}{6}\mathbf{Q}_{i+2,j}+\frac1{\Delta x}(\mathbf{Q}_{i+1,j}-\mathbf{Q}_{i,j})x\\
		&+\frac1{\Delta x^2}\left(\frac12\mathbf{Q}_{i,j}-\mathbf{Q}_{i+1,j}+\frac12\mathbf{Q}_{i+2,j}\right)x^2,\\
		p^{r5}_3(x)&=c_0+\frac1{\Delta x}c_1x+\frac1{\Delta x^2}c_2x^2+\frac1{\Delta x^3}c_3x^3+\frac1{\Delta x^4}c_4x^4.\\
	\end{aligned}
\end{equation}
where $c_0,\cdots, c_5$ are the coefficients, which have the form
$$
\begin{aligned}
	c_0 &= \frac{1}{60}(2\mathbf{Q}_{i-2,j}-13\mathbf{Q}_{i-1,j}+47\mathbf{Q}_{i,j}+27\mathbf{Q}_{i+1,j}-3\mathbf{Q}_{i+2,j}),\\
	c_1&=\frac1{12}(\mathbf{Q}_{i-1,j}-15\mathbf{Q}_{i,j}+15\mathbf{Q}_{i+1,j}-\mathbf{Q}_{i+2,j}),\\
	c_2&=\frac18(-\mathbf{Q}_{i-2,j}+6\mathbf{Q}_{i-1,j}-8\mathbf{Q}_{i,j}+2\mathbf{Q}_{i+1,j}+\mathbf{Q}_{i+2,j}),\\
	c_3&=\frac16(-\mathbf{Q}_{i-1,j}+3\mathbf{Q}_{i,j}-3\mathbf{Q}_{i+1,j}+\mathbf{Q}_{i+2,j}),\\
	c_4&=\frac1{24}(\mathbf{Q}_{i-2,j}-4\mathbf{Q}_{i-1,j}+6\mathbf{Q}_{i,j}-4\mathbf{Q}_{i+1,j}+\mathbf{Q}_{i+2,j}).
\end{aligned}
$$
Considering the weights of the stencils, we can rewrite $p^{r5}_3(x)$ as
\begin{equation}
	p^{r5}_3(x)=d_3\left(\frac{1}{d_3}p^{r5}_3(x)-\sum^2_{k=0}\frac{d_k}{d_3}p^{r3}_k(x)\right)+\sum^2_{k=0}d_kp^{r3}_k(x),\ d_3\neq0.
\end{equation}
where $d_k,\ k=0,1,2,3$ are defined as linear weights and have the following values
\begin{equation}
	d_3=d_{H},\quad d_0=(1-d_H)(1-d_L)/2,\quad d_1=(1-d_H)d_L,\quad d_2=d_0
\end{equation}
where $d_H\in[0.85,0.95]$ and $d_L\in[0.85,0.95]$, and here we choose $d_H=d_L=0.85$.

To deal with the discontinuity, the WENO-Z type \cite{borges2008improved} non-linear weights are used as
\begin{equation}
	\omega_k=d_k\left(1+\left(\frac{\tau_Z}{\beta_k+\varepsilon}\right)^2\right).
\end{equation}
where the global smoothness indicator $\tau_Z$ is defined as
\begin{equation}\label{eq37}
	\tau_Z=\frac13\left(|\beta_3-\beta_0|+|\beta_3-\beta_1|+|\beta_3-\beta_2|\right).
\end{equation} 
The local smoothness indicators are defined as
\begin{equation}\label{eq38}
	\beta_k=\sum^{q_k}_{q=1}\Delta x^{2q-1}\int^{x_{i+1/2,j}}_{x_{i-1/2,j}}\left(\frac{{\rm d}^q}{{\rm d}x^q}p_k(x)\right)^2{\rm d}x.
\end{equation}
where $q_k$ is the order of $p_k(x)$. Combined with the normalized weights $\overline{\omega}_k=\omega_k/(\sum^3_0\omega_k)$, the final form of the reconstruction polynomial is 
\begin{equation}
	P(x)=\overline{\omega}_3\left(\frac{1}{d_3}p^{r5}_3(x)-\sum^2_{k=0}\frac{d_k}{d_3}p^{r3}_k(x)\right)+\sum^2_{k=0}\overline{\omega}_kp^{r3}_k(x),\ d_3\neq0.
\end{equation}
After the normal reconstruction is completed to obtain the line value of the interface, the value of each Gaussian point at the interface can be obtained similarly after the tangential reconstruction. More details on WENO-AO reconstruction can be found in Ref \cite{ji2019high}.

\subsection{Discontinuity feedback factor}
The idea of DF is to start from the interface reconstruction values at the current time step and then predict the cell that the discontinuities will enter at the next time step, which is first proposed in \cite{ji2021gradient}, here we use the improved form \cite{yue2022high} and denote $\alpha_{ij}\in[0,1]$ as DF at a cell ${\rm\Omega}_{ij}$
\begin{equation}
	\alpha_{ij}=\prod_{p=1}^4\prod_{m=1}^2\alpha_{p,m}.
\end{equation}
where $\alpha_{p,m}$ is the DF corresponding to the Gaussian point $\mathbf{x}_{p,m}$, which has the form
\begin{equation}
	\begin{aligned}
		&\alpha_{p,m}=\frac1{1+D^2},\\
		&D=\frac{\left|p^l-p^r\right|}{p^l}+\frac{\left|p^l-p^r\right|}{p^r}+\left({\rm Ma}^l_n-{\rm Ma}^r_n\right)^2+\left({\rm Ma}^l_t-{\rm Ma}^r_t\right)^2.
	\end{aligned}
\end{equation}
where $p^k, k=l,r$ denote the left and right pressure of the Gaussian point, ${\rm Ma}_n^l, {\rm Ma}_t^l$ denote the left Mach number defined by the normal and tangential velocity.

\begin{figure}[htbp]
	\centering
	\includegraphics[width=1\textwidth]{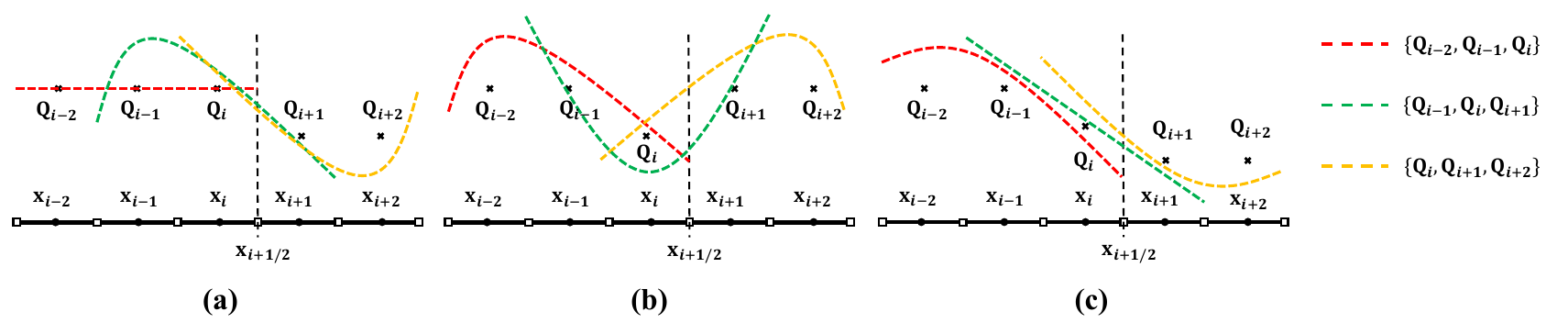}
	\caption{Three possible distributions of variables $\mathbf{Q}_{i-2}\cdots\mathbf{Q}_{i+2}$. (a) Smooth sub-stencils exist, and the WENO reconstruction will automatically approximated to $\{\mathbf{Q}_{i-2},\mathbf{Q}_{i-1},\mathbf{Q}_{i}\}$ by weights. (b-c) Each sub-stencil has a discontinuity, and the WENO reconstruction can only select the relatively smooth sub-stencil by weights, which means the effectiveness of the reconstruction polynomial for each sub-stencil is reduced. } 
	\label{fig:1}       
\end{figure} 

\subsection{Hybrid reconstruction}
When there is a discontinuity, the smootnest sub-stencil is automatically approximated by the weights, but if both the cell ${\rm \Omega}_{i,j}$ and its neighboring cells ${\rm \Omega}_{i-1,j},{\rm \Omega}_{i+1,j}$ have discontinuities, the reconstruction polynomials corresponding to each sub-stencil have no effect (see Fig.~\ref{fig:1}). For example, if there are discontinuities only in cell ${\rm \Omega}_{i,j}$ and cell ${\rm \Omega}_{i+1,j}$, the WENO-Z reconstruction will approximate the reconstruction value to the smoothest sub-stencil $\{{\rm \Omega}_{i-2,j},{\rm \Omega}_{i-1,j},{\rm \Omega}_{i,j}\}$. If there are discontinuities in the cell ${\rm \Omega}_{i-1,j},{\rm \Omega}_{i,j},{\rm \Omega}_{i+1,j}$, any sub-stencil and large-stencil will not be effective, leading to a decrease in the robustness of the algorithm. To improve the robustness of the algorithm, we include a DF threshold $\alpha_{thres}$ to determine the existence of a discontinuity within the cell, which means that if $\alpha_{i,j}<\alpha_{thres}$, there is a discontinuity in the cell.

\begin{figure}[htbp]
	\centering
	\includegraphics[width=1\textwidth]{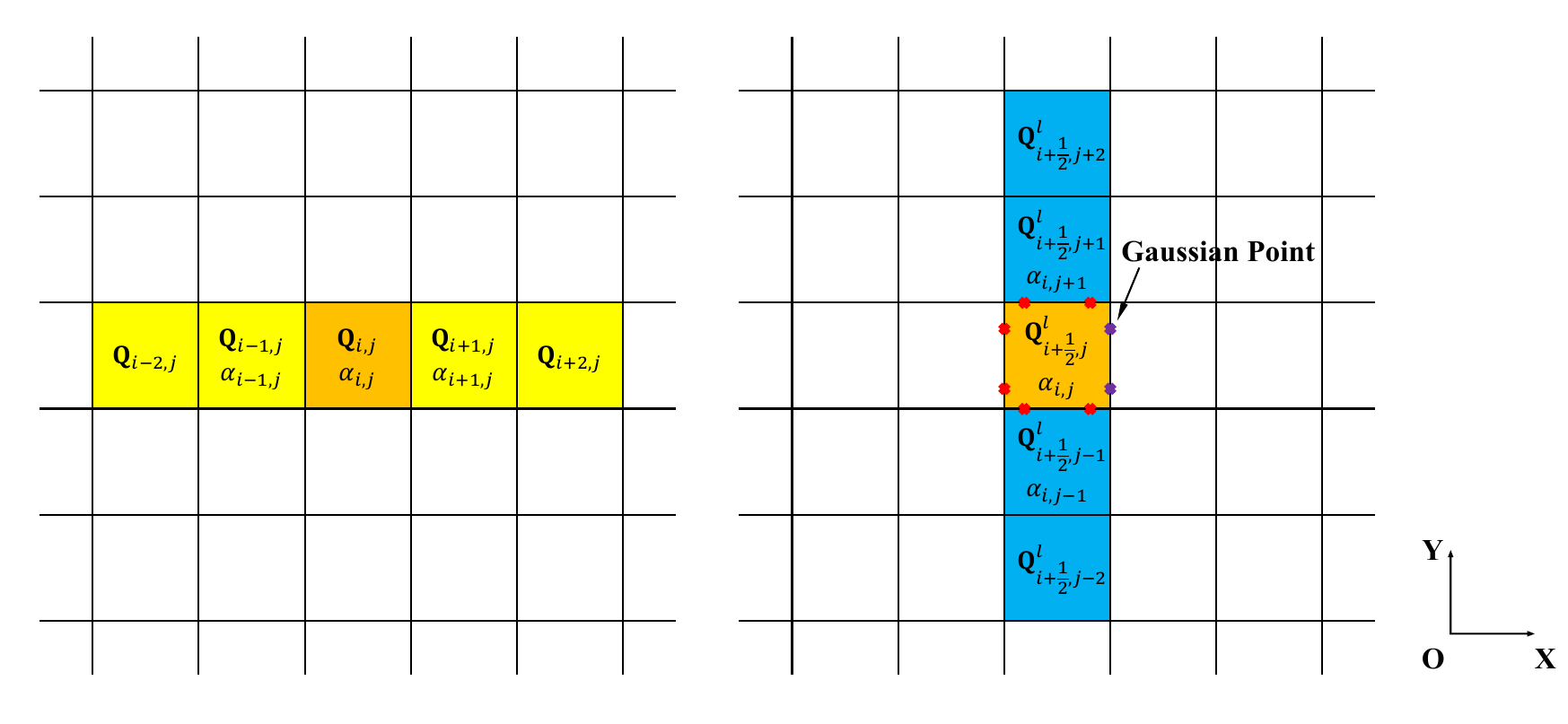}
	\caption{x-direction reconstruction at the left side of the interface ${\rm \Gamma}_{i+1/2,j}$. Normal velocity and tangential velocity(from left to right).}
	\label{fig:22}       
\end{figure}

\paragraph{Remark 1} The DF corresponds to the strength of the discontinuity. As long as the DF is 1, the flow is smooth. However, the numerical solution itself is not strictly continuous, which means that some cells have a DF of less than 1, but there is no physical discontinuity in the cell. To deal with this drawback, we set a threshold $\alpha_{thres}$ for the DF, which is used to determine whether there is a discontinuity in the cell. If the threshold is close to 0, it means more discontinuity, and if the threshold is close to 1, it means higher resolution. To improve the robustness and maintain the high resolution of the algorithm, $\alpha_{thres}=0.5$ is used in this paper.

For the cell ${\rm \Omega}_{i,j}$, if the following condition is satisfied
\begin{equation}\label{eq34}
	{\rm max}\{\alpha_{i-1,j},\ \alpha_{i,j},\ \alpha_{i+1,j}\}<\alpha_{thres}
\end{equation}
we change the reconstruction method from WENO-AO to the quadratic polynomial $p^{r3}_1(x)$, convert $p^{r3}_1(x)$ to zero-mean form and add DF, we have the reconstruction form
\begin{equation}
	\begin{aligned}
		p^{r3}_{DF}(x)&=\mathbf{Q}_{i,j}\\
		&+\alpha_{i,j}\left[\frac{1}{\Delta x}(\mathbf{Q}_{i+1,j}-\mathbf{Q}_{i,j})(x-x_0)+\frac{1}{\Delta x^2}(\frac12\mathbf{Q}_{i-1,j}-\mathbf{Q}_{i,j}+\frac12\mathbf{Q}_{i+1,j})(x^2-x_1)\right].
	\end{aligned}
\end{equation}
where $x_0=-\frac12\Delta x$ and $x_1=\frac13\Delta x^2$. Subject to satisfying Eq. \ref{eq34}, if there is a discontinuity in cell ${\rm \Omega}_{ij}$, $\alpha_{i,j}$ will be approximated to zero, and the reconstruction value of the interface ${\rm\Gamma}_{i+1/2,j}$ will be approximated to $\mathbf{Q}_{i,j}$.

Since the DF is defined on the cell, the value of the DF of the cell on the same side as the Gaussian point of the interface is used for the tangential reconstruction(see Fig.~\ref{fig:22}).

\begin{algorithm}[htbp]\label{Alo1}
	\SetAlgoLined 
	\caption{High-order finite volume scheme with the hybrid reconstruction algorithm}
	\KwIn{Cell average values $\mathbf{W}^n$, DF value $\alpha^n$ at $t_n$ and $\alpha_{threshold}$}
	\KwOut{Cell average values $\mathbf{W}^{n+1}$ and DF value $\alpha^{n+1}$ at $t_{n+1}$}
	\BlankLine
	{calculate} time increment $\Delta t\leftarrow(\mathbf{W}^n,CFL)$\;
	\While{not at end of stages}{
		{update} boundary condition\;
		\ForEach{cell ${\rm \Omega}_{ij}$}{
			calculate characteristic variables $\mathbf{Q}_{ij}\leftarrow \mathbf{W}_{ij}$\;
			\tcp{as an example, interface $(i+1/2,j)$ left side reconstruction is shown.}
			\If{without condition {\rm min}\{{\rm max}\{$\alpha_{i-1,j},\alpha_{i,j},\alpha_{i+1,j}$\}$,\alpha_{threshold}$\}}{cycle termination;}
			\tcp{normal reconstruction}
			\If{x-direction}{
				\eIf{{\rm max}\{$\alpha_{i-1,j},\alpha_{i,j},\alpha_{i+1,j}$\}$<\alpha_{threshold}$ }{
					using 3rd polynomial reconstruction with DF\;
					calculate  $(\mathbf{Q}^l_{i+1/2,j},(\mathbf{Q}_x)^l_{i+1/2,j})\leftarrow(\mathbf{Q}_{i-1,j},\mathbf{Q}_{i,j},\mathbf{Q}_{i+1,j},\alpha_{i,j})$
				}{
					using 5th WENO-AO reconstruction\;
					calculate $(\mathbf{Q}^l_{i+1/2,j},(\mathbf{Q}_x)^l_{i+1/2,j})\leftarrow(\mathbf{Q}_{i-2,j},\mathbf{Q}_{i-1,j},\mathbf{Q}_{i,j},\mathbf{Q}_{i+1,j},\mathbf{Q}_{i+2,j})$
				}
			}
			\If{y-direction}{similar to x-direction}
			\tcp{tangential reconstruction}
			\If{x-direction}{
				\ForEach{Gaussian point $\mathbf{x}_m$ at the interface $(i+1/2,j)$}{
					\tcp{here the index $i+1/2$ is omitted}
					\eIf{{\rm max}\{$\alpha_{i,j-1},\alpha_{i,j},\alpha_{i,j+1}$\}$<\alpha_{threshold}$ }{
						using 3rd polynomial reconstruction with DF\;
						calculate  $(\mathbf{Q}^l_{j,m},(\mathbf{Q}_x)^l_{j,m},(\mathbf{Q}_y)^l_{j,m})\leftarrow(\mathbf{Q}^l_{j-1},\mathbf{Q}^l_{j},\mathbf{Q}^l_{j+1},(\mathbf{Q}_x)^l_{j},\alpha_{i,j})$
					}{
						using 5th WENO-AO reconstruction\;
						calculate $(\mathbf{Q}^l_{j,m},(\mathbf{Q}_x)^l_{j,m},(\mathbf{Q}_y)^l_{j,m})\leftarrow(\mathbf{Q}^l_{j-2},\mathbf{Q}^l_{j-1},\mathbf{Q}^l_{j},\mathbf{Q}^l_{j+1},\mathbf{Q}^l_{j+2},(\mathbf{Q}_x)^l_{j})$
					}
				}
			}
			\If{y-direction}{similar to x-direction}
			calculate $(\mathbf{W}^l_{j,m},(\mathbf{W}_x)^l_{j,m},(\mathbf{W}_y)^l_{j,m})\leftarrow(\mathbf{Q}^l_{j,m},(\mathbf{Q}_x)^l_{j,m},(\mathbf{Q}_y)^l_{j,m})$\;
		}
		\tcp{similar to interface (i+1/2,j) right side reconstruction.}
		\ForEach{interface}{
			calculate interface equilibrium state values and first-order derivatives $(\mathbf{W}^c_{i+1/2,j},(\mathbf{W}_x)^c_{i+1/2,j})\leftarrow(\mathbf{W}^{l,r}_{i+1/2,j},(\mathbf{W}_x)^{l,r}_{i+1/2,j})$,\\
			\ForEach{Gaussian point $\mathbf{x}_{m}$ at the interface $(i+1/2,j)$}{
				calculate $(\mathbf{W}^c_{j,m},(\mathbf{W}_x)^c_{j,m},(\mathbf{W}_y)^c_{j,m})$\;
				calculate flux $\mathbf{F}_{j,m}$ and first-order derivatives $\partial_t\mathbf{F}_{j,m}$.
			}
		}
		\ForEach{cell ${\rm \Omega}_{ij}$}{
			update conservative variables $\mathbf{W}^*_{i,j}$ at $t_*$.
		}
	}
	{update} DF $\alpha^{n+1}$ at $t_{n+1}$\;
	{update} current time $t_{n+1}=t_n+\Delta t$\;
	{return}
\end{algorithm}

\subsection{Reconstruction of equilibrium state for gas-kinetic scheme}
For the non-equilibrium state, the reconstruction can be used to obtain the equilibrium state $g^c,g^c_x$ and $g^c_y$ directly by $g^k,g^k_\mathbf{x},k=l,r$, and here a kinetic-based weighting method is used
\begin{equation}
	\begin{aligned}
		\int\bm{\psi}g^c{\rm d}\Xi&=\mathbf{W}^c=\int_{u>0}\bm{\psi}g^l{\rm d}\Xi+\int_{u<0}\bm{\psi}g^r{\rm d}\Xi,\\
		\int\bm{\psi}g^c_\mathbf{x}{\rm d}\Xi&=\mathbf{W}^c_\mathbf{x}=\int_{u>0}\bm{\psi}g^l_\mathbf{x}{\rm d}\Xi+\int_{u<0}\bm{\psi}g^r_\mathbf{x}{\rm d}\Xi,\\
	\end{aligned}
\end{equation}
In this way, all components in Eq. (\ref{eq24}) can be determined. The procedure of the robustness enhanced high-order GKS algorithm is shown in Algorithm \ref{Alo1}.

\section{Numerical Examples}
\label{sec:5}
In this section, high-order finite volume schemes using the hybrid reconstruction method and two flux solvers are used to simulate several classical Riemann problems. The first scheme is the 2rd-order GKS solver, using the S2O4 time discretization with the non-linear limiter \cite{zhao2023direct}, which can limit the possible discontinuity of the flux function in time; the other scheme is the Lax-Friedrichs(L-F) solver using the strong stability preserving RK(SSP-RK) time discretization, which has the positivity preserving property. With the above two flux solvers, we could compare the resolution and robustness of the hybrid reconstruction with the WENO-AO reconstruction.
\subsection{Accuracy validations}
For the Euler equation, the smooth sin-wave propagation is used for the accuracy evaluation \cite{ji2018family}. In these cases, both the physical viscosity and the collision time are set to zero. The initial condition is given by
\begin{equation}
	\rho(x)=1+0.2{\rm sin}(\pi x),\ U(x)=1.0,\ p(x)=1.0,\ x\in[0,2].
\end{equation}
and the exact solution with periodic boundary conditions is 
\begin{equation}
	\rho(x,t)=1.0+0.2{\rm sin}(\pi(x-t)),\ U(x,t)=1.0, \ p(x,t)=1.0.
\end{equation}
Here we choose a time step of $\Delta t=0.25\Delta x$ to test the spatial and temporal accuracy together, which corresponds to $CFL\approx 0.5$. The numerical solutions from the hybrid scheme are obtained after a periodic of propagation at $t=2.0$, and compared with the exact result at $t=2.0$. Based on the $L^1, L^2$ and $L^{\infty}$ errors, the orders of the hybrid scheme using two solvers at $t=2.0$ are presented in Table~\ref{tab:1}-\ref{tab:LF1D}.

Extending the accuracy test to the 2-D case, the initial condition of the 2-D sin-wave propagation is 
\begin{equation}
	\begin{aligned}
	 	&\rho(x,y)=1.0+0.2{\rm sin}(\pi x){\rm sin}(\pi y),\\
		&U(x,y) = 1.0,\ V(x,y)=1.0,\ p(x,y)=1.0.
	\end{aligned}
\end{equation}
with the exact solution
\begin{equation}
	\begin{aligned}
		&\rho(x,y,t)=1.0+0.2{\rm sin}(\pi(x-t)){\rm sin}(\pi(y-t)),\\
		&U(x,y,t)=1.0,\ V(x,y,t)=1.0,\ p(x,y,t)=1.0.
	\end{aligned}
\end{equation}
The computational domain is $[-1,1]\times[-1,1]$ and the periodic boundary conditions is used. A CFL number $CFL=0.1$ is used for accuracy test cases. The results are shown in Table~\ref{tab:2}-\ref{tab:LF2D}. The results show that both the 1-D case and the 2-D case achieve their theoretical numerical accuracy.

\begin{table}[htbp]
\caption{Accuracy test for the 1-D sin-wave propagation using S2O4 GKS solver with hybrid reconstruction}
\label{tab:1}       
\begin{tabular*}{\textwidth}{@{\extracolsep{\fill}}ccccccc}
\hline\noalign{\smallskip}
Mesh number & $L^1$ error & order  & $L^2$ error & order & $L^\infty$ error & order\\
\noalign{\smallskip}\hline\noalign{\smallskip}
20 & 3.086057e-05 & & 3.429008e-05 & & 5.035811e-05 & \\
40 & 9.734233e-07 & 4.99 & 1.078578e-06 &4.99 & 1.593003e-06 &4.98 \\
80 & 3.045127e-08 & 5.00 & 3.375302e-08 & 5.00& 4.992874e-08 & 5.00\\
160 & 9.518025e-10 & 5.00 & 1.054893e-09 & 5.00& 1.561347e-09 & 5.00\\
\noalign{\smallskip}\hline
\end{tabular*}
\end{table}

\begin{table}[htbp]
	\caption{Accuracy test for the 1-D sin-wave propagation using SSP-RK3 L-F solver with hybrid reconstruction}
	\label{tab:LF1D}       
	\begin{tabular*}{\textwidth}{@{\extracolsep{\fill}}ccccccc}
		\hline\noalign{\smallskip}
		Mesh number & $L^1$ error & order  & $L^2$ error & order & $L^\infty$ error & order\\
		\noalign{\smallskip}\hline\noalign{\smallskip}
		20 & 8.721602e-05 & & 9.691440e-05 & & 1.435365e-04  & \\
		40 & 2.778818e-06 & 4.97 & 3.094034e-06 &4.87 & 4.627061e-06 &4.96 \\
		80 & 8.825108e-08 & 4.98 & 9.823738e-08 & 4.98& 1.470886e-07 & 4.98\\
		160 & 2.882101e-09  & 4.94 & 3.207652e-09 & 4.94& 4.793567e-09 & 4.94\\
		\noalign{\smallskip}\hline
	\end{tabular*}
\end{table}

\begin{table}[htbp]
	\caption{Accuracy test for the 2-D sin-wave propagation using S2O4 GKS solver with hybrid reconstruction}
	\label{tab:2}       
	\begin{tabular*}{\textwidth}{@{\extracolsep{\fill}}ccccccc}
		\hline\noalign{\smallskip}
		Mesh number & $L^1$ error & order  & $L^2$ error & order & $L^\infty$ error & order\\
		\noalign{\smallskip}\hline\noalign{\smallskip}
		20$\times$20 & 5.659621e-05 &  & 6.287560e-05 & & 9.141807e-05 & \\
		40$\times$ 40& 1.799309e-06 & 4.98 & 1.994645e-06 & 4.98 & 2.917438e-06 & 4.97\\
		80$\times$80 & 5.778633e-08 & 4.96 & 6.408196e-08 & 4.96& 9.346639e-08 & 4.96\\
		160$\times$160 & 1.968348e-09 & 4.88 & 2.186488e-09  & 4.87& 3.155101e-09  & 4.89\\
		\noalign{\smallskip}\hline
	\end{tabular*}
\end{table}

\begin{table}[htbp]
	\caption{Accuracy test for the 2-D sin-wave propagation using SSP-RK3 L-F solver with hybrid reconstruction}
	\label{tab:LF2D}       
	\begin{tabular*}{\textwidth}{@{\extracolsep{\fill}}ccccccc}
		\hline\noalign{\smallskip}
		Mesh number & $L^1$ error & order  & $L^2$ error & order & $L^\infty$ error & order\\
		\noalign{\smallskip}\hline\noalign{\smallskip}
		20$\times$20 & 1.744519e-04 &  & 1.939279e-04 & & 4.420067e-04 & \\
		40$\times$ 40& 5.643100e-06  & 4.95 & 6.282656e-06 & 4.95 & 9.380738e-06 & 4.92\\
		80$\times$80 & 1.872392e-07 & 4.91 & 2.083394e-07  & 4.91& 3.108420e-07 & 4.94\\
		160$\times$160 & 7.075964e-09 & 4.73 & 7.871208e-09 & 4.73& 1.163168e-08 & 4.73\\
		\noalign{\smallskip}\hline
	\end{tabular*}
\end{table}

\subsection{Resolution validations for non-linear waves}
To validate the resolution of different reconstruction methods under the two solvers, here three Riemann problems are used to study the complicated non-linear wave structures.
\subsubsection{Shu-Osher problem}
\begin{figure}[htbp]
	\centering
	\subfigure[S2O4+GKS, Global]{
		\includegraphics[width=0.48\textwidth]{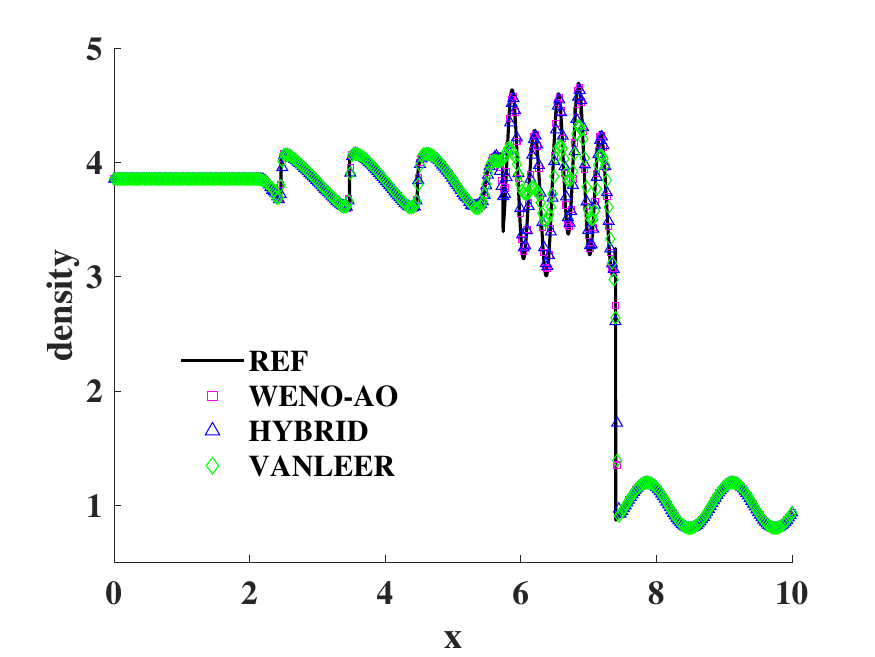}
	}
	\subfigure[S2O4+GKS, Local]{
		\includegraphics[width=0.48\textwidth]{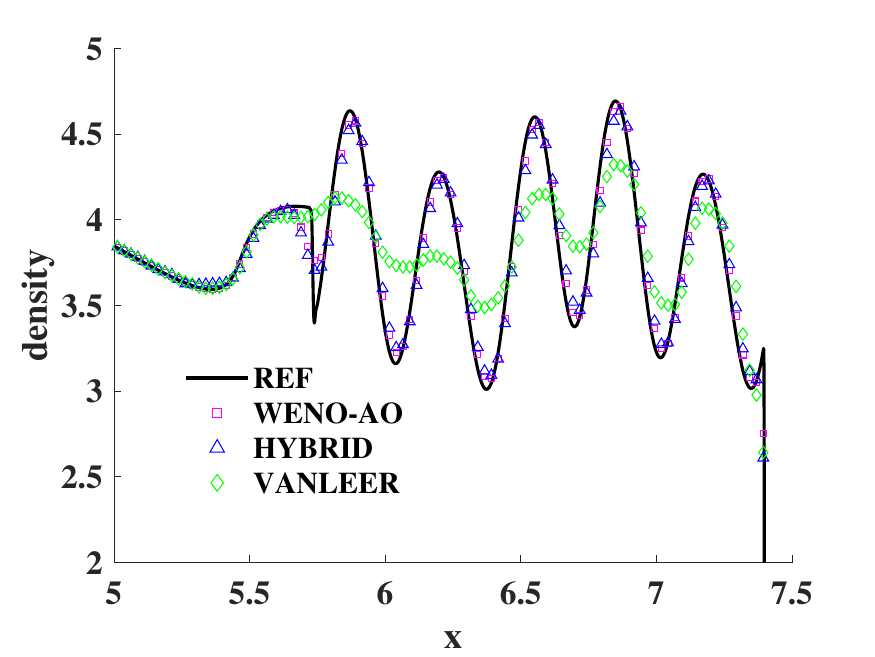}
	}
	\subfigure[SSP-RK3+L-F, Global]{
		\includegraphics[width=0.48\textwidth]{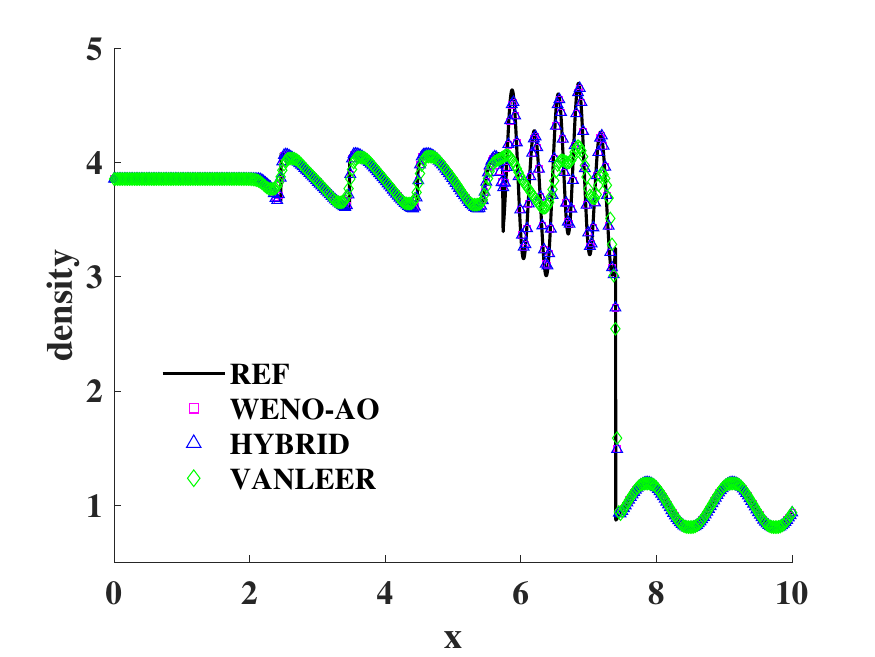}
	}
	\subfigure[SSP-RK3+L-F, Local]{
		\includegraphics[width=0.48\textwidth]{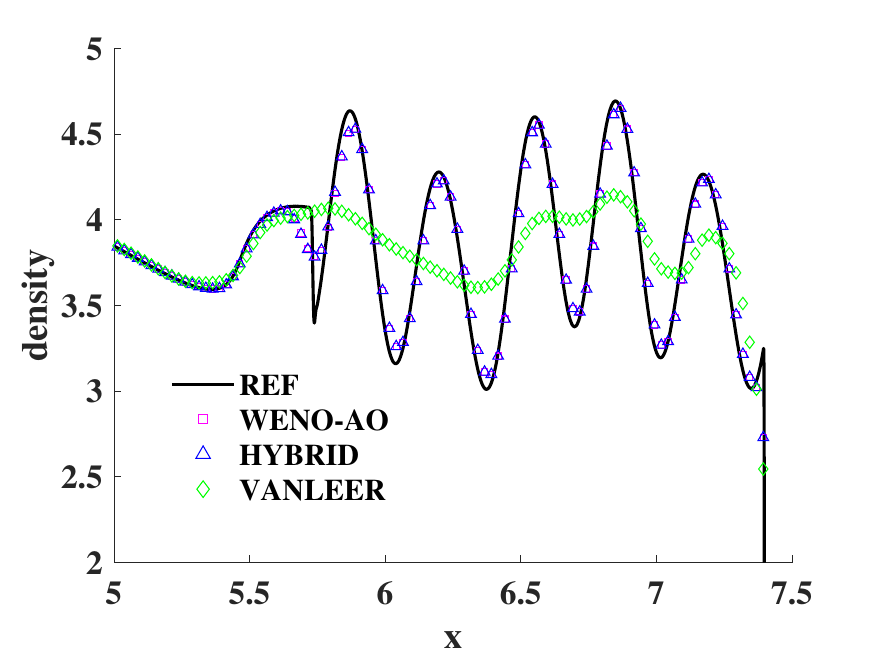}
	}
	\caption{Shu-Osher problem: the density distributions and local enlargement at $t=1.8$ with a cell size $\Delta x=1/40$. (a-b) Different reconstruction methods using the S2O4 GKS solver. (c-d) Different reconstruction methods using the SSP-RK3 L-F solver. The reference solution is obtained by the 1-D fifth-order WENO-AO GKS with 10000 meshes.}
	\label{Shu-Osher}
\end{figure}

The 1-D Shu-Osher problem \cite{shu1989efficient} can test the performance of capturing high frequency wave, and the initial condition is given as follows
\begin{equation}
	(\rho,u,p)=
	\begin{cases}
		(3.857134, 2.629369, 10.33333),\quad &x\in[0.0,1.0],\\
		(1.0+0.2{\rm sin}(5x),0,1),\quad &x\in[1.0,10.0].
	\end{cases}
\end{equation}
The computational domain is $[0.0,10.0]$ with a mesh size of $\Delta x=1/40$. A CFL number $CFL=0.5$ is used for the all following cases. Numerical results of the density distributions under the two solvers using different reconstruction methods are shown in Fig.~\ref{Shu-Osher} at $t=1.8$. For the same mesh size, both the spatial accuracy of the WENO-AO and the hybrid reconstruction are higher than that of the van Leer reconstruction, with the former performing better than the latter in capturing extreme values. Furthermore, the performance of the algorithm for capturing extreme values using the hybrid reconstruction is slightly lower than using the WENO-AO reconstruction but still maintains a good resolution.

\subsubsection{Interaction of planar shocks}
\begin{figure}[htbp]
	\centering
	\subfigure[S2O4+GKS, WENO-AO]{
		\includegraphics[width=0.48\textwidth]{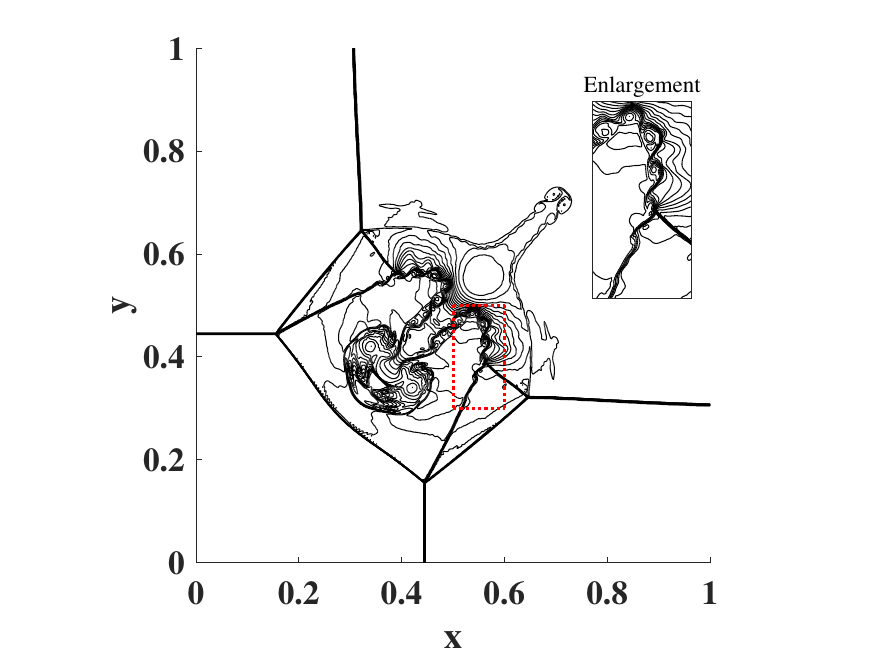}
	}
	\subfigure[S2O4+GKS, Hybrid]{
		\includegraphics[width=0.48\textwidth]{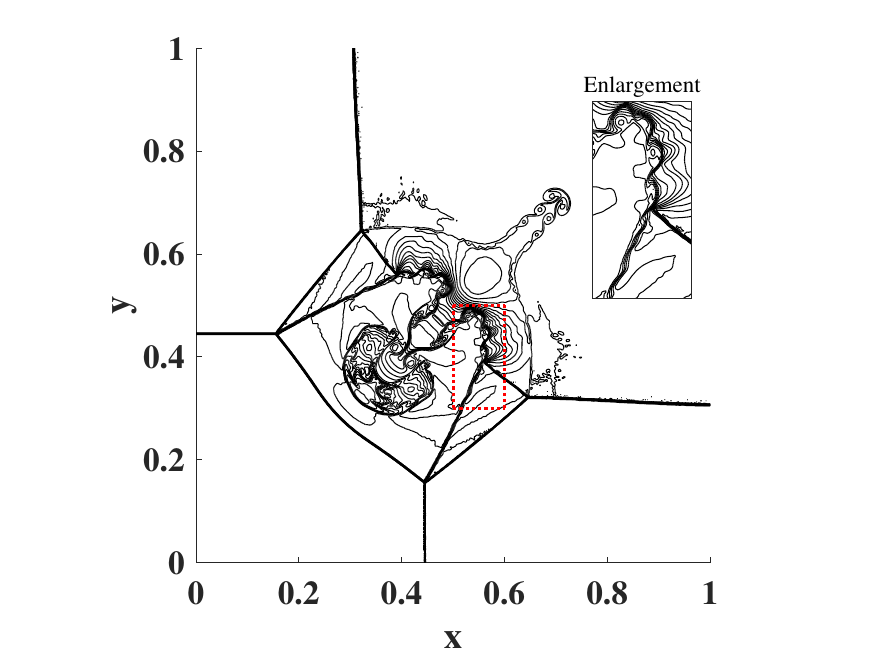}
	}
	\subfigure[SSP-RK3+L-F, WENO-AO]{
		\includegraphics[width=0.48\textwidth]{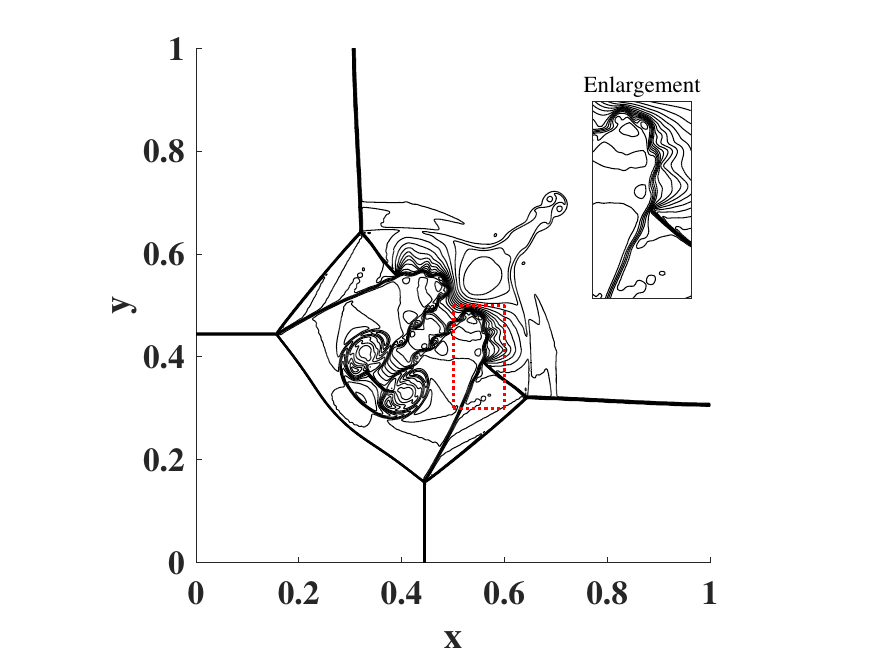}
	}
	\subfigure[SSP-RK3+L-F, Hybrid]{
		\includegraphics[width=0.48\textwidth]{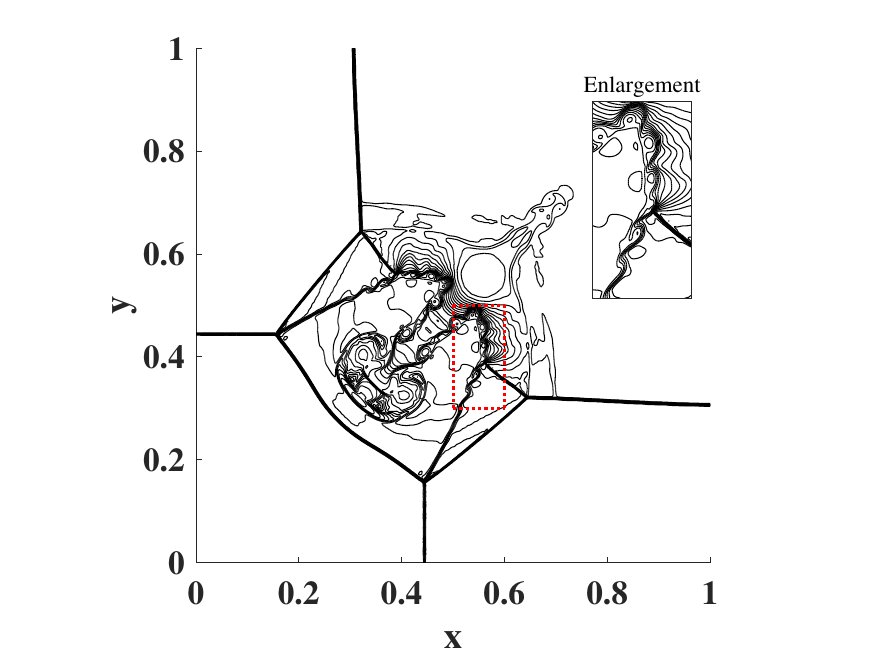}
	}
	\caption{Interaction of planar shocks problem: the density distributions and local enlargement at $t=0.6$ with $500\times 500$ meshes. (a-b) Different reconstruction methods using S2O4 GKS solver. (c-d) Different reconstruction methods using SSP-RK3 L-F solver.}
	\label{planar shocks}
\end{figure}

Configuration 3 in Ref \cite{1998Solution} includes the shokc-shock interaction and the shock-vortex interaction. The initial condition is given as follows
\begin{equation}
	(\rho,u,v,p)=
	\begin{cases}
		(0.138,1.206,1.206,0.129),\quad &x<0.7,y<0.7,\\
		(0.5323,0,1.206,0.3),&x\ge 0.7,y<0.7,\\
		(1.5,0,0,1.5),&x\ge 0.7,y\ge 0.7,\\
		(0.5323,1.206,0,0.3),&x<0.7,y\ge 0.7.
	\end{cases}
\end{equation}
Using the same 25 density contours, the numerical solutions using different reconstruction methods at $t=0.6$ are shown in Fig~\ref{planar shocks}. Both reconstruction methods can capture the shock sharply, compared to WENO-AO, the hybrid method has a higher resolution in the shear layer, i.e. in the enlargement region in Fig~\ref{planar shocks}, and shows stronger instabilities of the vortex sheets in the enlargement region, suggesting that the new algorithm can maintain the high resolution.

\subsubsection{Interaction of planar contact discontinuity}
\begin{figure}[htbp]
	\centering
	\subfigure[S2O4+GKS, WENO-AO]{
		\includegraphics[width=0.48\textwidth]{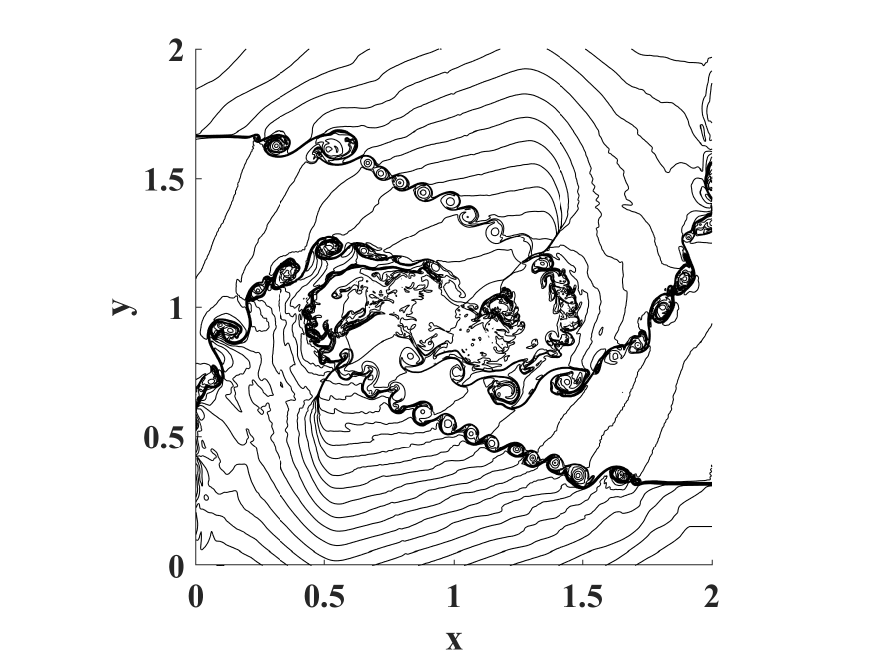}
	}
	\subfigure[S2O4+GKS, Hybrid]{
		\includegraphics[width=0.48\textwidth]{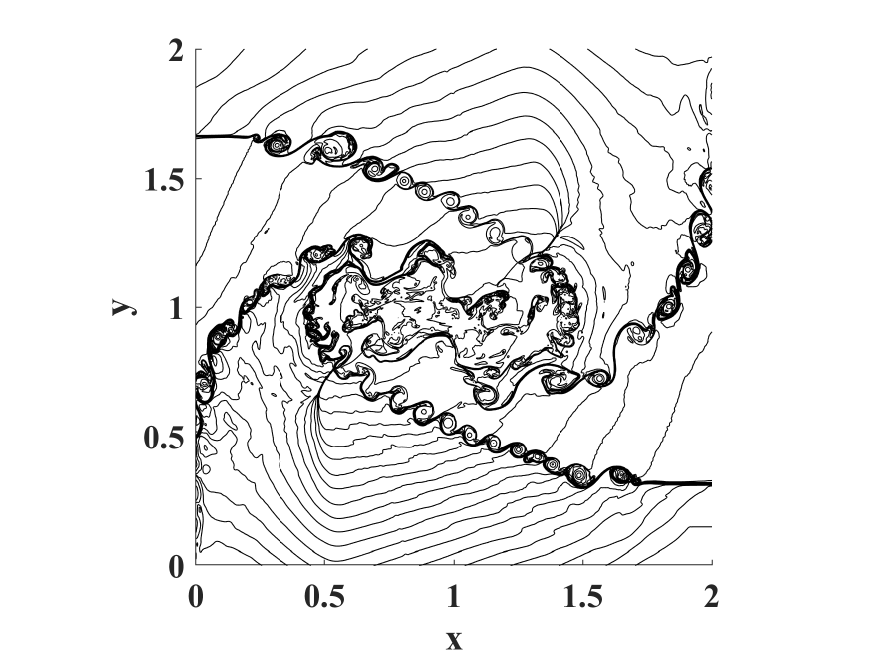}
	}
	\subfigure[SSP-RK3+L-F, WENO-AO]{
		\includegraphics[width=0.48\textwidth]{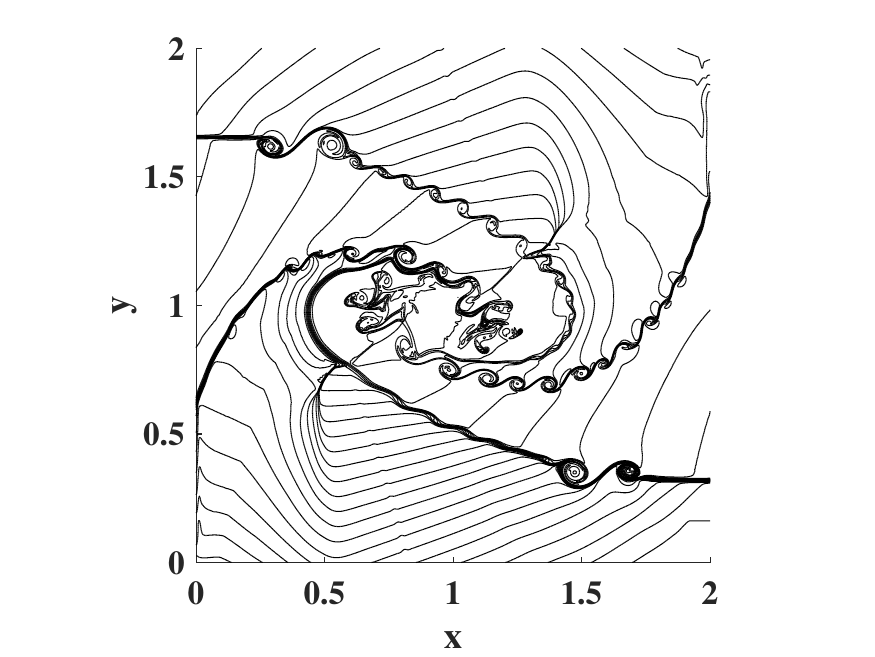}
	}
	\subfigure[SSP-RK3+L-F, Hybrid]{
		\includegraphics[width=0.48\textwidth]{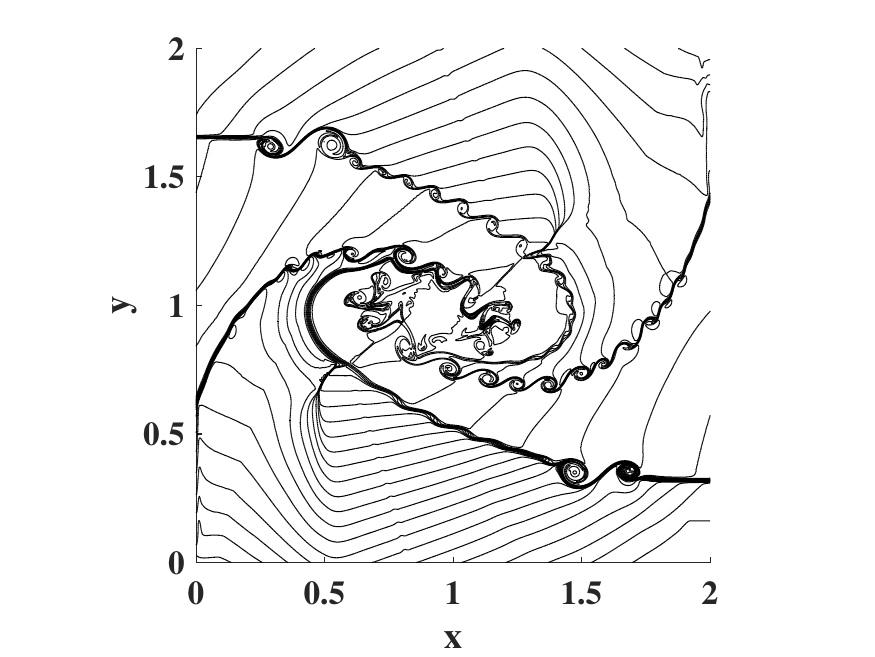}
	}
	\caption{Interaction of planar contact discontinuous: the density distributions at $t=1.6$ with $800\times 800$ meshes. (a-b) Different reconstruction methods using S2O4 GKS solver. (c-d) Different reconstruction methods using SSP-RK3 L-F solver.}
	\label{planar contact}
\end{figure}

Configuration 6 in \cite{1998Solution} involves planar contact discontinuous interactions. The initial condition in a computational domain $[0,2]\times[0,2]$ is given by
\begin{equation}
	(\rho,u,v,p)=
	\begin{cases}
		(1.0,-0.75,0.5,1.0),\quad &x<1.0,y<1.0,\\
		(3.0,-0.75,-0.5,1.0),&x\ge 1.0,y<1.0,\\
		(1.0,0.75,-0.5,1.0),&x\ge 1.0,y\ge 1.0,\\
		(2.0,0.75,0.5,1.0),&x<1.0,y\ge 1.0.
	\end{cases}
\end{equation}
As shown in Fig \ref{planar contact}, the flow structure is complicated, very close results are obtained by the WENO-AO/hybrid reconstruction, both the reconstruction methods show high resolution, such as many small vortexes.

\subsection{Robustness validations for non-linear waves}
Three Riemann problems are used to evaluate the robustness of the algorithm.
\subsubsection{123 problem}
\begin{figure}[htbp]
	\centering
	\subfigure[WENO-AO]{
		\includegraphics[width=0.31\textwidth]{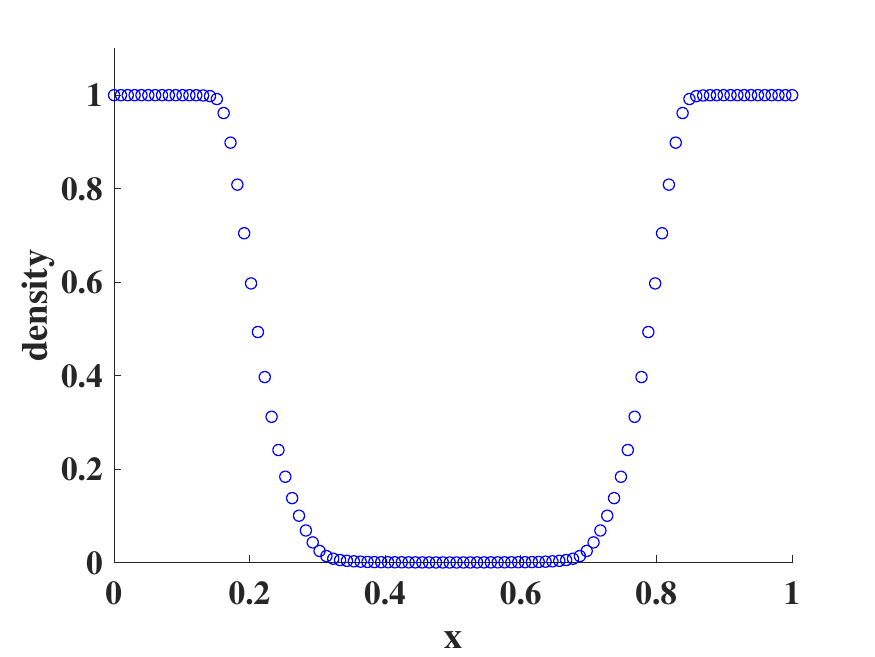}
	}
	\subfigure[Hybrid]{
		\includegraphics[width=0.31\textwidth]{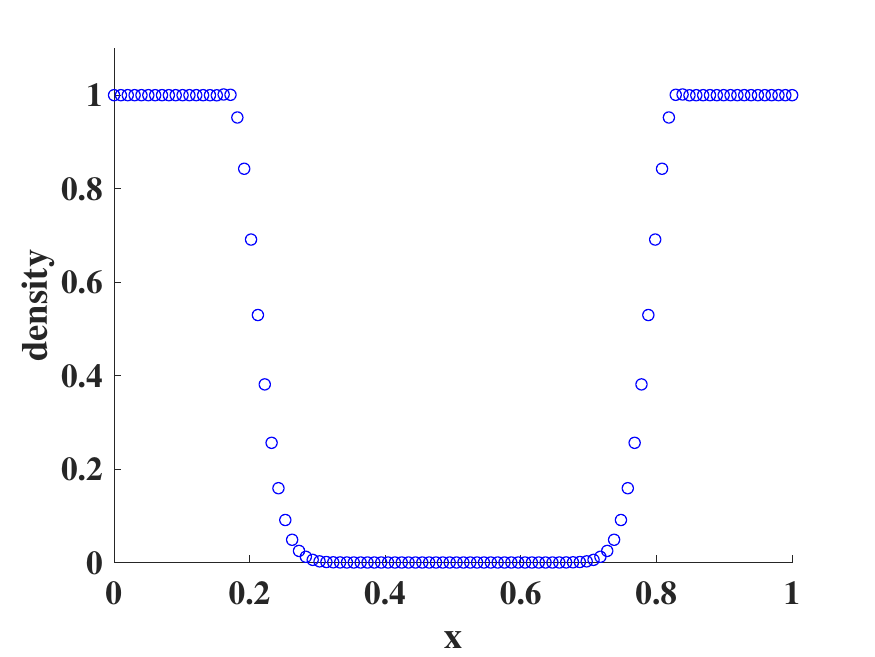}
	}
	\subfigure[van Leer]{
		\includegraphics[width=0.31\textwidth]{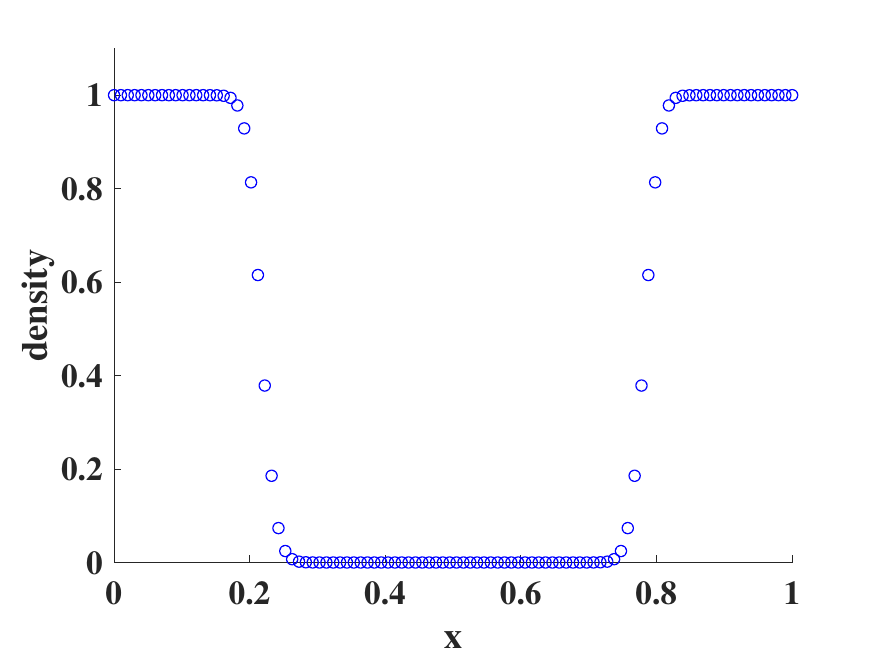}
	}
	\caption{123 problem: density distributions using different reconstruction methods at the corresponding maximum Mach number under GKS solver. Left: Ma= 6.0 with original WENO-AO reconstruction. Middle: Ma=31.5 with hybrid reconstruction. Right: Ma=31.5 with van Leer reconstruction.}
	\label{fig:4}
\end{figure}

\begin{figure}[htbp]
	\centering
	\subfigure[WENO-AO]{
		\includegraphics[width=0.31\textwidth]{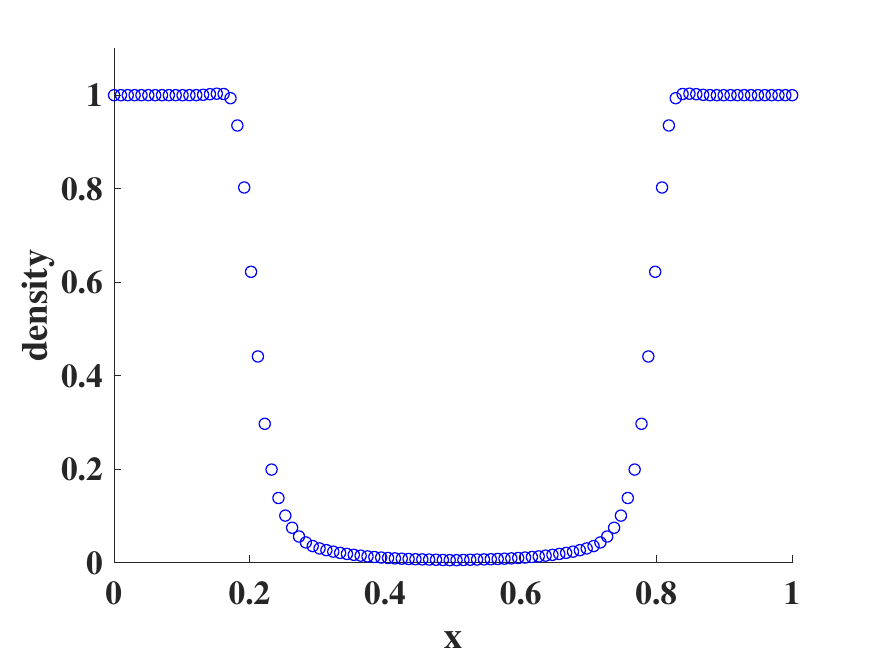}
	}
	\subfigure[Hybrid]{
		\includegraphics[width=0.31\textwidth]{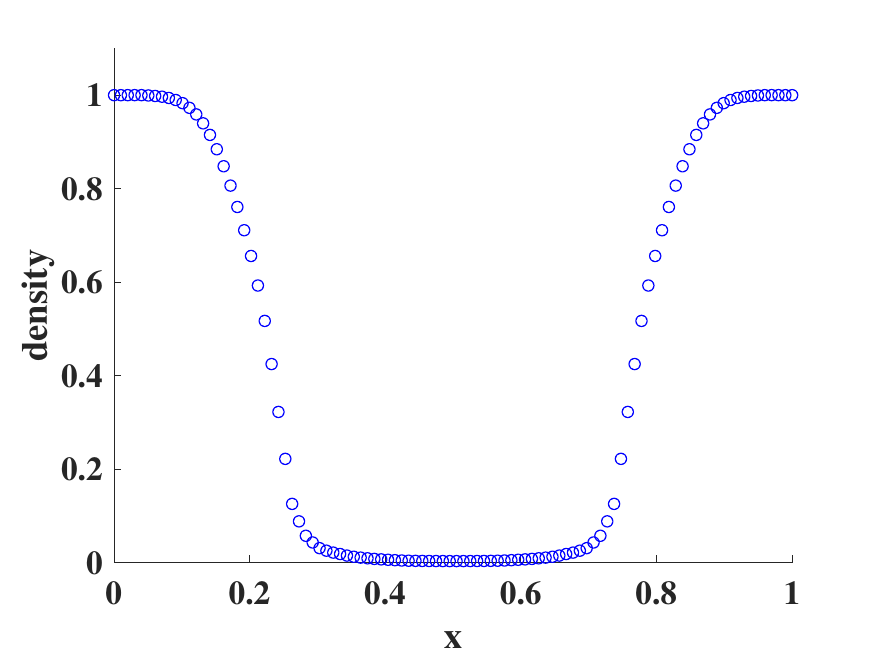}
	}
	\subfigure[van Leer]{
		\includegraphics[width=0.31\textwidth]{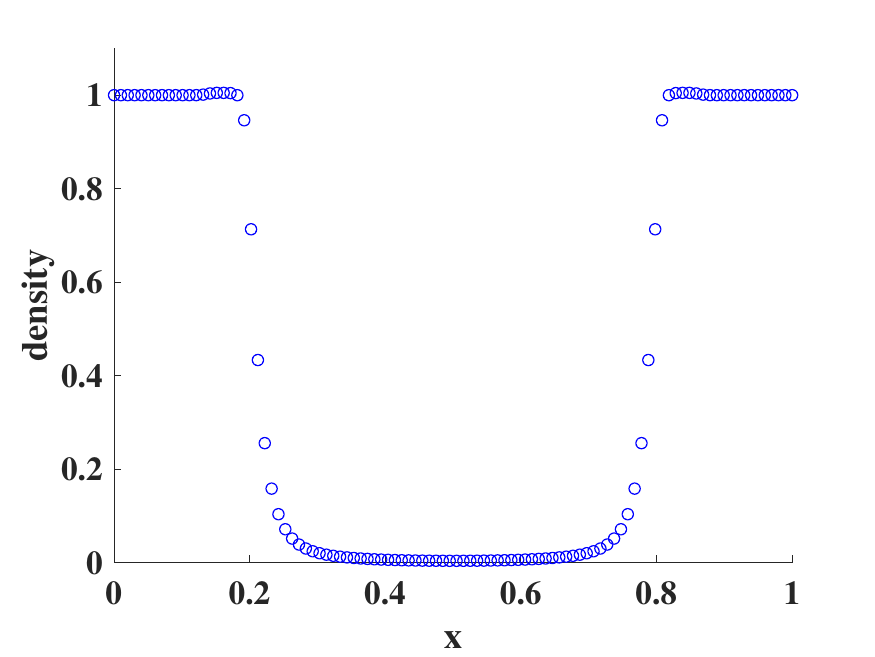}
	}
	\caption{123 problem: density distributions using different reconstruction methods at the corresponding maximum Mach number under L-F solver. Left: Ma= 7.9 with original WENO-AO reconstruction. Middle: Ma=84.1 with hybrid reconstruction. Right: Ma=119.5 with van Leer reconstruction.}
	\label{fig:5}
\end{figure}

The 123 problem \cite{liska2003comparison} yields two strong rarefaction waves propagating in both directions, with the density and pressure in the center region approaching zero. The initial condition is given as follows
\begin{equation}
	(\rho,u,p)=\begin{cases}
		(1.0,-2.0,p_0),\quad &x\in[0,0.5)\\
		(1.0,2.0,p_0),&x\in[0.5,1.0]
	\end{cases}
\end{equation}
In this case, the reference Mach number can be obtained by adjusting the value of $p_0$, which has the form
\begin{equation}
	\mathbf{Ma}=\frac{u}{c}=\frac{u}{\sqrt{{\gamma p_0}/{\rho}}}
\end{equation}
where $\gamma$ is the specific heat ratio, and the maximum Mach number that the algorithm can achieve by testing is used to evaluate the robustness of the algorithm. Free boundary conditions and 100 meshes are used, and the solutions using different solvers with hybrid reconstruction at $t=0.14$ are shown in Fig~\ref{fig:4}-\ref{fig:5} and Table~\ref{tab:3}. The results show that the maximum Mach number that can be computed by the algorithm is significantly higher when using the hybrid method than the WENO-AO under the GKS solver, and the hybrid reconstruction method achieves the same robustness as using the van Leer limiter, similarly under the L-F solver.
\begin{table}[htbp]
	\caption{123 problem: Maximum Mach number using different reconstruction methods}
	\label{tab:3}       
	\begin{tabular*}{\textwidth}{@{\extracolsep{\fill}}cccc}
		\hline\noalign{\smallskip}
		GKS solver & Maximum Mach number & LF solver & Maximum Mach number\\
		\noalign{\smallskip}\hline\noalign{\smallskip}
		WENO-AO & 6.0 & WENO-AO & 7.9\\
		Hybrid& 31.5 & Hybrid & 84.1\\
		van Leer & 31.5 & van Leer & 119.5\\
		\noalign{\smallskip}\hline
	\end{tabular*}
\end{table}

\subsubsection{Hurricane-like problem}
\begin{figure}[htbp]
	\centering
	\subfigure[WENO-AO]{
		\includegraphics[width=0.51\textwidth]{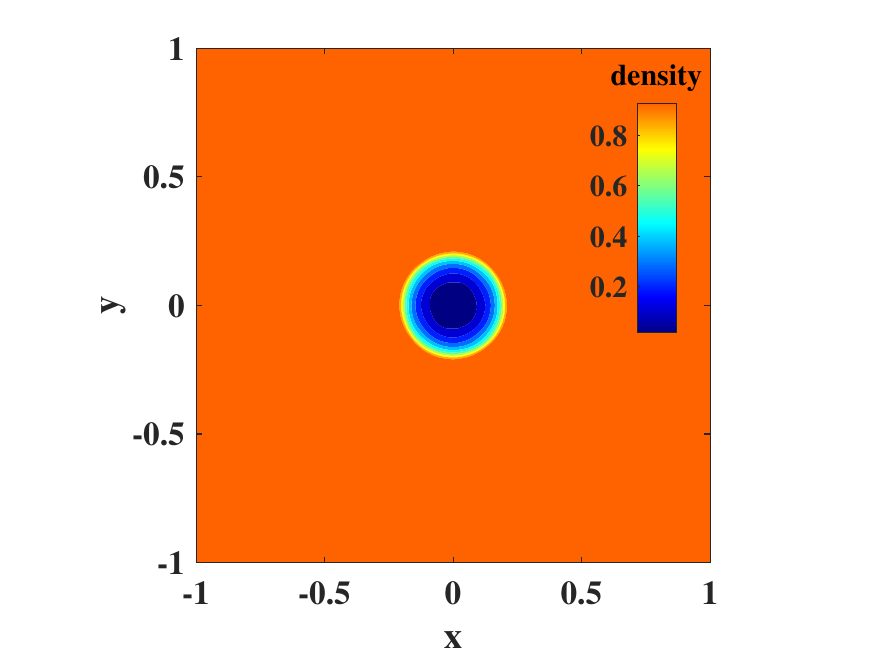}
	}\hspace{-10mm}
	\subfigure[Hybrid]{
		\includegraphics[width=0.51\textwidth]{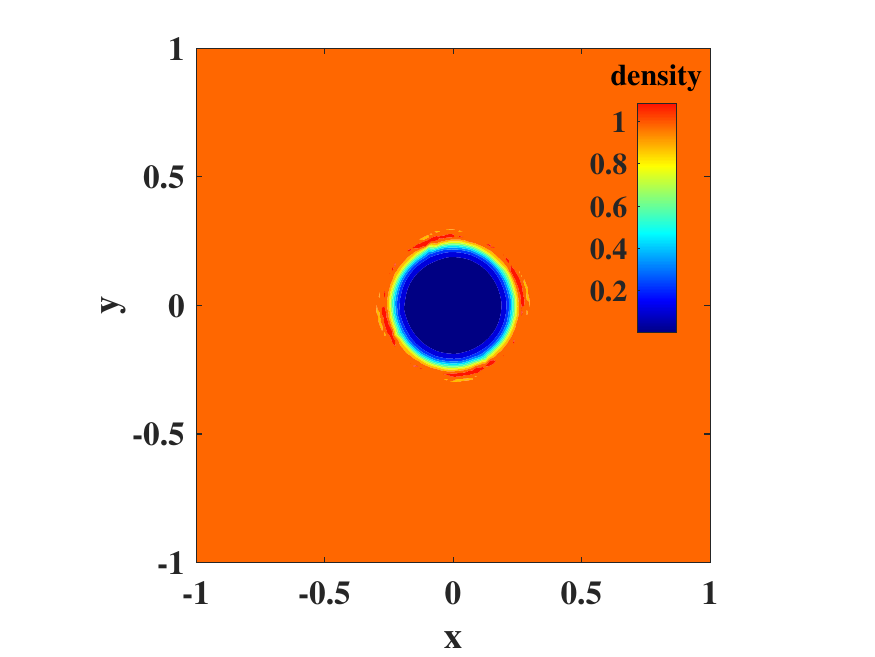}
	}
	\caption{Hurricane-like problem: density distributions using different reconstruction methods at the corresponding maximum Mach number under S2O4 GKS solver. Left: Ma= 2.0 with original WENO-AO reconstruction. Right: Ma=16.0 with hybrid reconstruction.}
	\label{fig:6}
\end{figure}

\begin{figure}[htbp]
	\centering
	\subfigure[WENO-AO]{
		\includegraphics[width=0.51\textwidth]{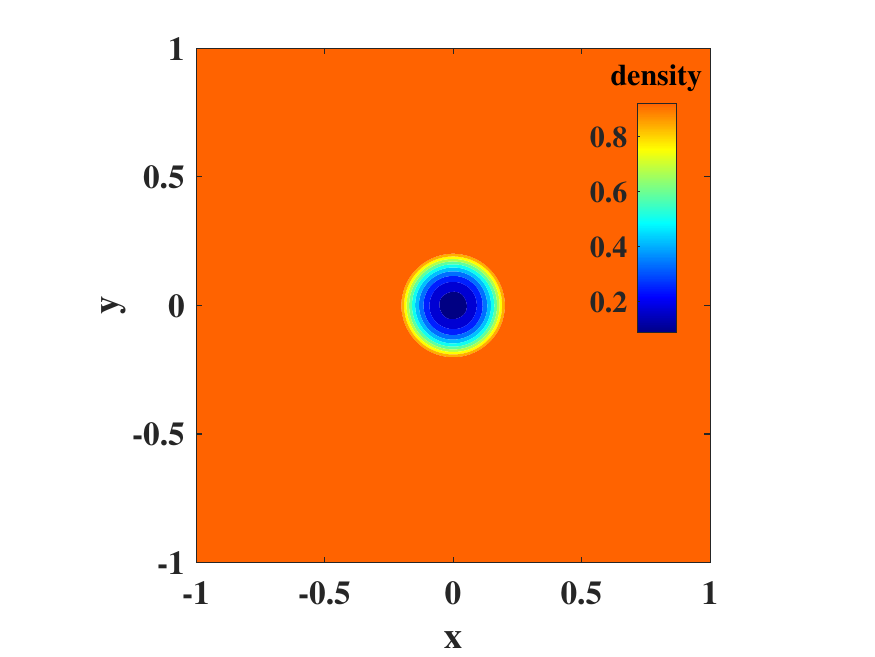}
	}\hspace{-10mm}
	\subfigure[Hybrid]{
		\includegraphics[width=0.51\textwidth]{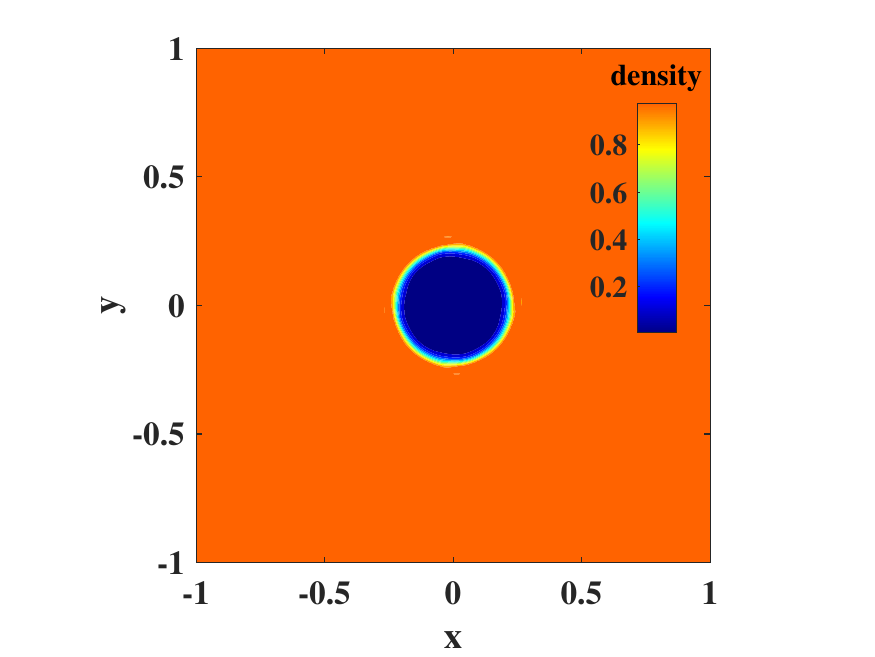}
	}
	\caption{Hurricane-like problem: density distributions using different reconstruction methods at the corresponding maximum Mach number under SSP-RK L-F solver. Left: Ma= 1.4 with original WENO-AO reconstruction. Right: Ma=7.0 with hybrid reconstruction.}
	\label{fig:7}
\end{figure}

The hurricane-like flow evolution \cite{zhang1996exact} has one-point vacuum in the center with rotational velocity field. The initial condition is given as follows
\begin{equation}
	(\rho,u,v,p)=(\rho_0,v_0{\rm sin}\ \theta,-v_0{\rm cos}\ \theta, A\rho_0^\gamma)
\end{equation}
where $\theta={\rm arctan}(y/x)$, $\gamma=1.4$, $\rho_0=1.0$, and $A=25$ is the initial entropy. The solutions are classified into three types according to the initial Mach number $\mathbf{Ma}=|v_0|/c_0$, where $c_0$ is the speed of sound. To evaluate the robustness, here we consider the high-speed rotation with $\mathbf{Ma}>\sqrt{2}$. In this case, the density goes to the vacuum faster and the fluid rotates strongly, which is a significant challenge for the robustness of the algorithm. For this case, we increase the Mach number by adjusting the value of $v_0$, and since no exact solution can be imposed on the boundary, the computational domain is $[-2,2]\times[-2,2]$, the mesh size is $\Delta x=\Delta y=1/100$, the non-reflection boundary is used. The solutions in the domain $[-1,1]\times[-1,1]$, which can ignore the effect of the boundary condition, are shown in Fig~\ref{fig:6}-\ref{fig:7} at the 50th time step.
\paragraph{Remark 2}In this case, when the initial velocity $v_0$ is increased, the rarefaction wave will propagate a longer distance in the same simulation time, and interact with the boundary conditions which leads to unphysical solutions. 
Furthermore, the strongest rarefaction wave is generated during the early simulation stage.
Thus, if the scheme can survive within the first 50 time steps, it will keep positivity in longer time simulation.\\

\begin{figure}[htbp]
	\centering
	\subfigure{
		\includegraphics[width=0.33\textwidth]{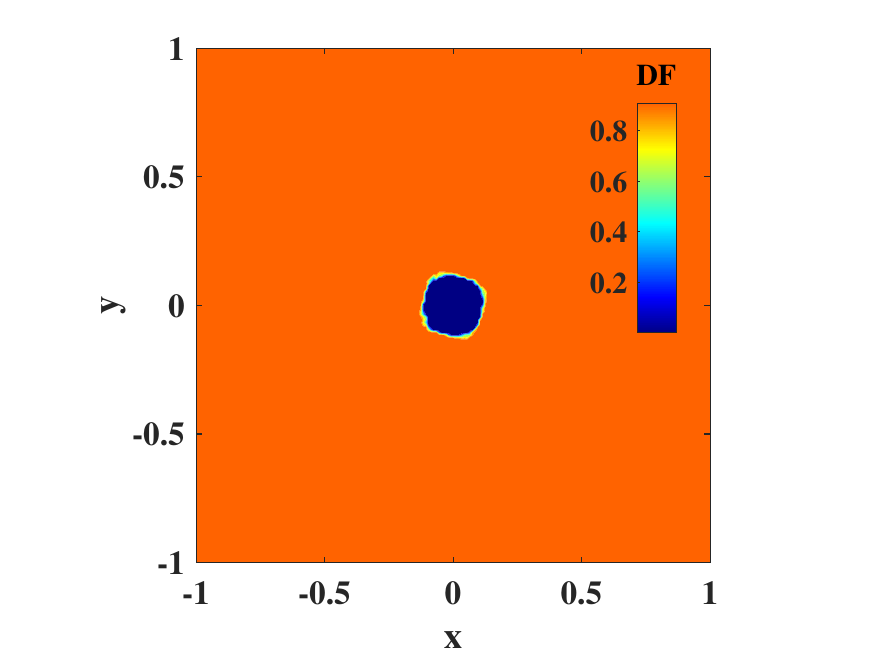}
	}\hspace{-5mm}
	\subfigure{
		\includegraphics[width=0.33\textwidth]{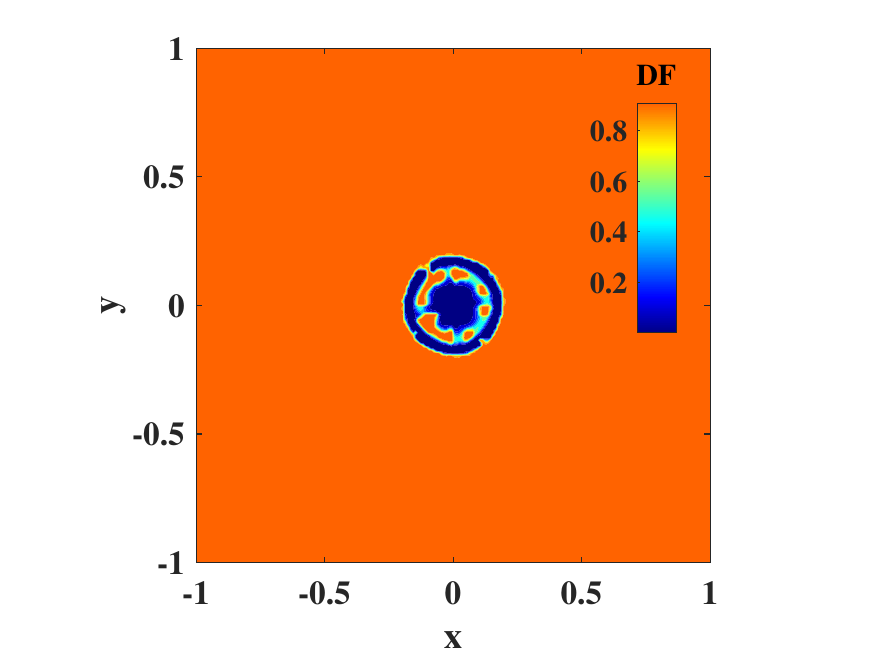}
	}\hspace{-5mm}
	\subfigure{
		\includegraphics[width=0.33\textwidth]{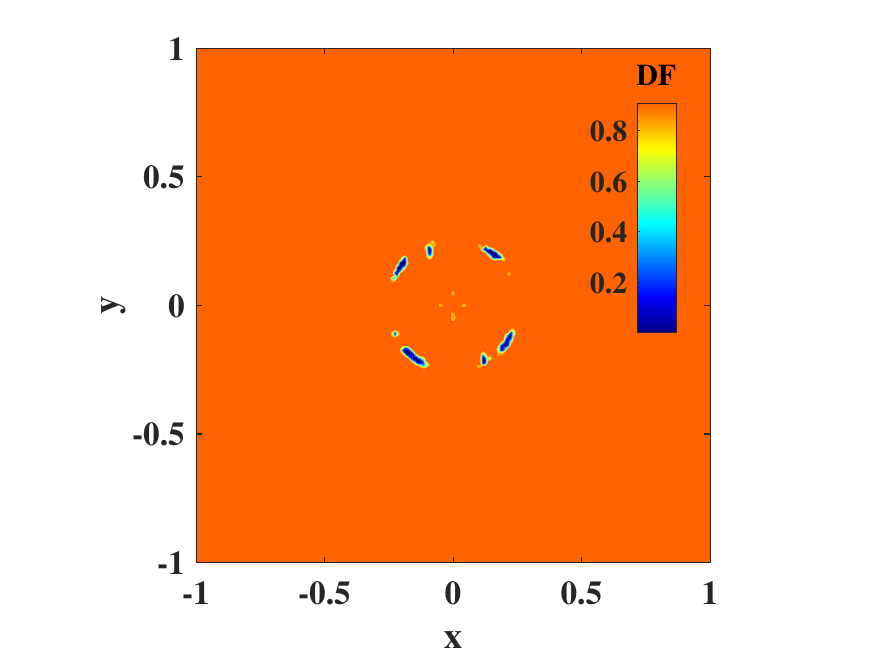}
	}
	\caption{DF distribution at step 15, 30, 45(from left to right) using S2O4 GKS solver.}
	\label{fig:8}
\end{figure}

 The time-dependent DF distributions at different steps are shown in Fig~\ref{fig:8}. At the early step, the slopes are modified near the center of the vacuum, where a rarefaction wave is immediately formed. As the rarefaction wave propagates outwards, the modified region of the DF extends outwards.
 The test shows that the current DF can improve the robustness of the hybrid under both the GKS and L-F solvers, as shown in Table~\ref{tab:4}.


\begin{table}[htbp]
	\caption{Hurricane-like problem: Maximum Mach number using different reconstruction methods}
	\label{tab:4}       
	\begin{tabular*}{\textwidth}{@{\extracolsep{\fill}}cccc}
		\hline\noalign{\smallskip}
		GKS solver & Maximum Mach number & LF solver & Maximum Mach number\\
		\noalign{\smallskip}\hline\noalign{\smallskip}
		WENO-AO & 2.0 & WENO-AO & 1.4\\
		Hybrid& 16.0 & Hybrid & 7.0\\
		\noalign{\smallskip}\hline
	\end{tabular*}
\end{table}

\subsubsection{Interaction of planar rarefaction waves}
Configuration 2 in Ref \cite{liska2003comparison} is about the interaction  of 2-D planar rarefaction waves, which will emerge a continuous transition from smooth flow to the presence of transonic shock. The initial condition is given as
\begin{equation}
	(\rho,u,v,p)=\begin{cases}
		(1,-0.6323,-0.6323,1.5), \quad &x<0.5,y<0.5\\
		(\rho,0.6323,-0.6323, A\rho^\gamma),\quad &x\ge 0.5,y<0.5\\
		(1.0.6323, 0.6323, 1.5),\quad &x\ge 0.5,y\ge 0.5\\
		(\rho,-0.6323, 0.6323, A\rho^\gamma),\quad &x<0.5,y\ge 0.5
	\end{cases}
\end{equation}
where $A$ is the initial entropy and $A=1.5$ in this case. For this case, the free boundary condition is used, the mesh size is $\Delta x=\Delta y=1/400$, and the computational domain is $[0,1]\times[0,1]$. Here we evaluate the robustness of the algorithm in the extreme case by decreasing the value of $p$. The reference Mach number is defined as 
\begin{equation}
	\mathbf{Ma}=\frac{\sqrt{U^2+V^2}}{\sqrt{\gamma A\rho^\gamma/\rho}}=\frac{\sqrt{U^2+V^2}}{\sqrt{\gamma A\rho^{\gamma-1}}}
\end{equation}

\begin{figure}[htbp]
	\centering
	\subfigure{
		\includegraphics[width=0.48\textwidth]{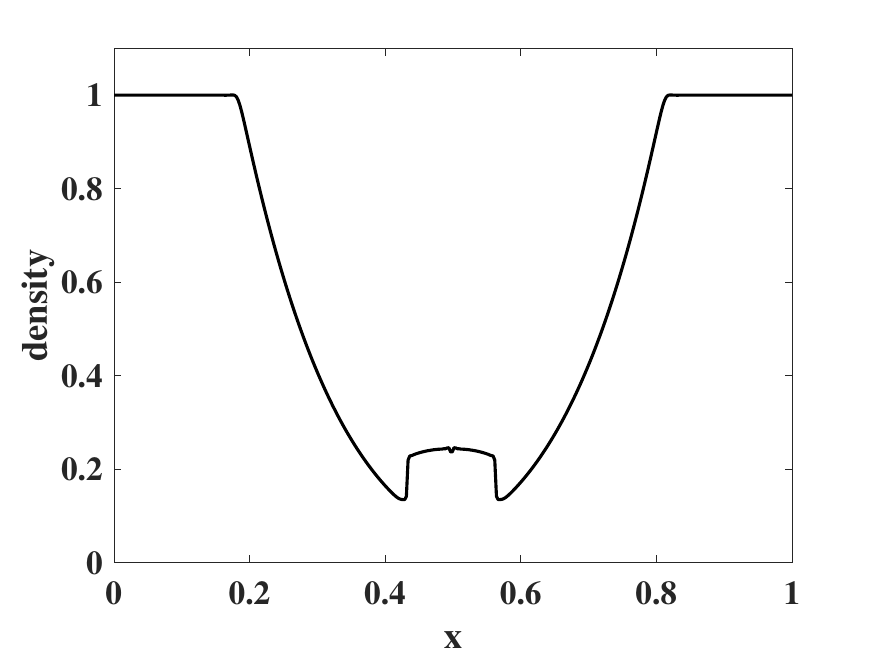}
	}
	\subfigure{
		\includegraphics[width=0.48\textwidth]{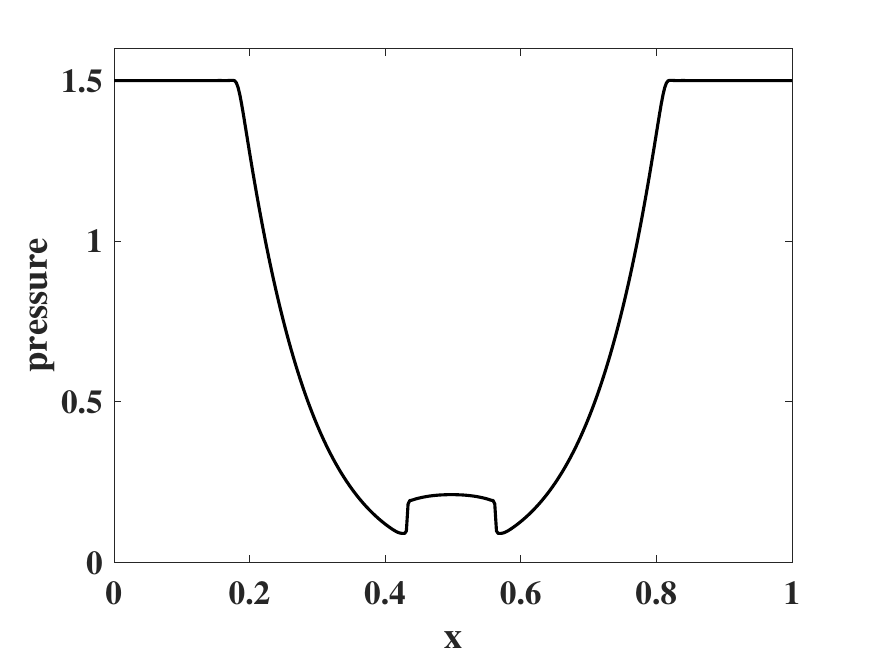}
	}
	\caption{The density and pressure distributions along the diagonal line. Using GKS solver with $p=0.4$.}
	\label{fig:24}
\end{figure}

\begin{figure}[htbp]
	\centering
	\subfigure[WENO-AO]{
		\includegraphics[width=0.51\textwidth]{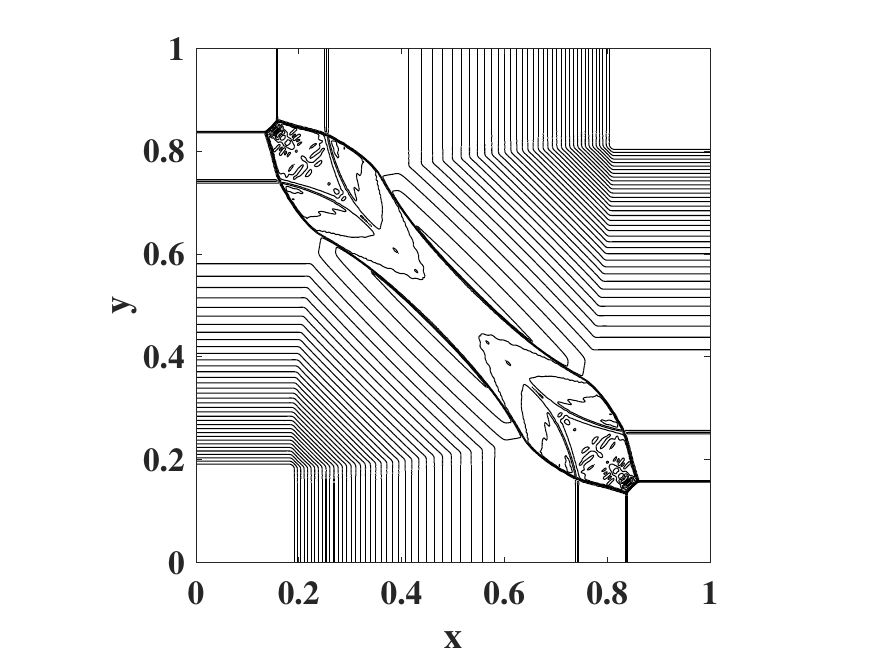}
	}\hspace{-10mm}
	\subfigure[Hybrid]{
		\includegraphics[width=0.51\textwidth]{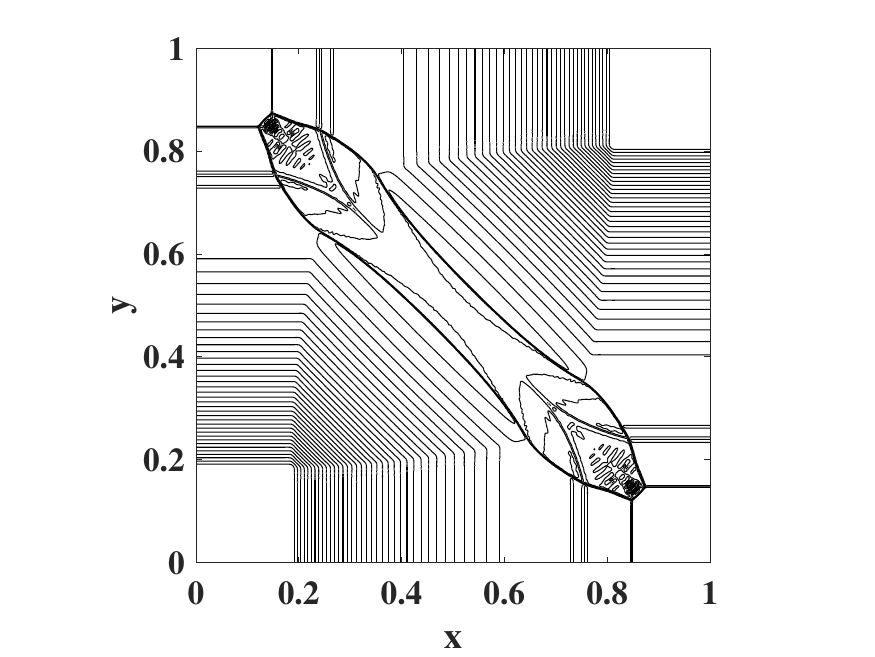}
	}
	\caption{Interaction of planar rarefaction waves: density distributions using different reconstruction methods at the corresponding maximum Mach number under S2O4 GKS solver. Left: Ma= 1.8 with original WENO-AO reconstruction. Right: Ma=2.1 with hybrid reconstruction.}
	\label{fig:25}
\end{figure}

\begin{figure}[htbp]
	\centering
	\subfigure[WENO-AO]{
		\includegraphics[width=0.51\textwidth]{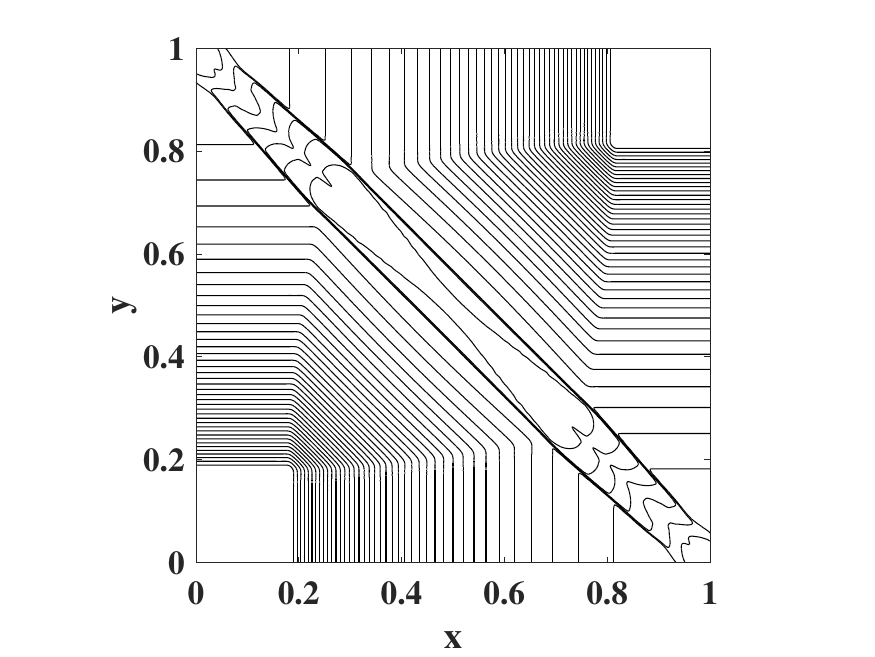}
	}\hspace{-10mm}
	\subfigure[Hybrid]{
		\includegraphics[width=0.51\textwidth]{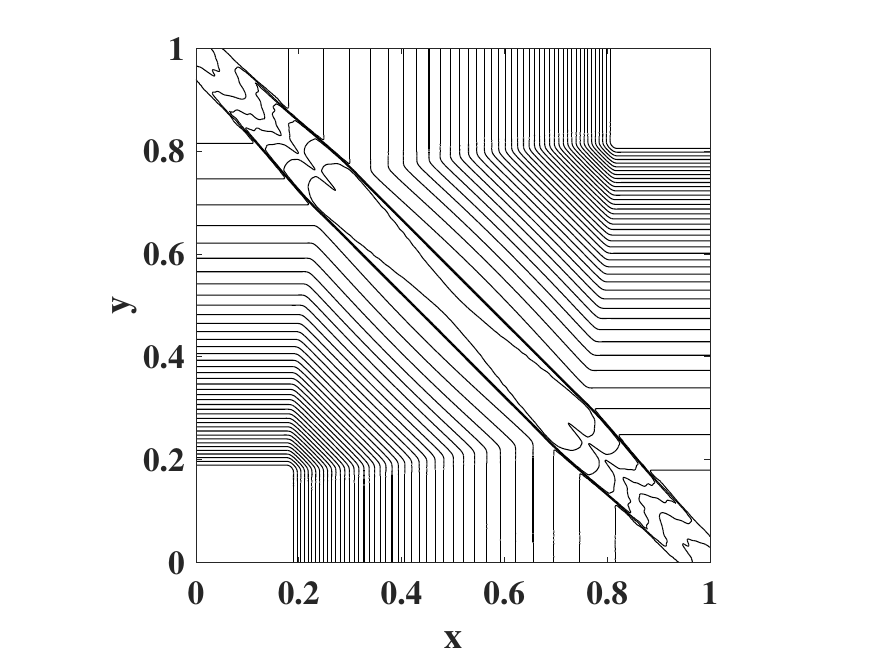}
	}
	\caption{Interaction of planar rarefaction waves: density distributions using different reconstruction methods at the corresponding maximum Mach number under SSP-RK L-F solver. Left: Ma= 18.1 with original WENO-AO reconstruction. Right: Ma=24.8 with hybrid reconstruction.}
	\label{fig:26}
\end{figure}

\begin{figure}[htbp]
	\centering
	\subfigure{
		\includegraphics[width=0.33\textwidth]{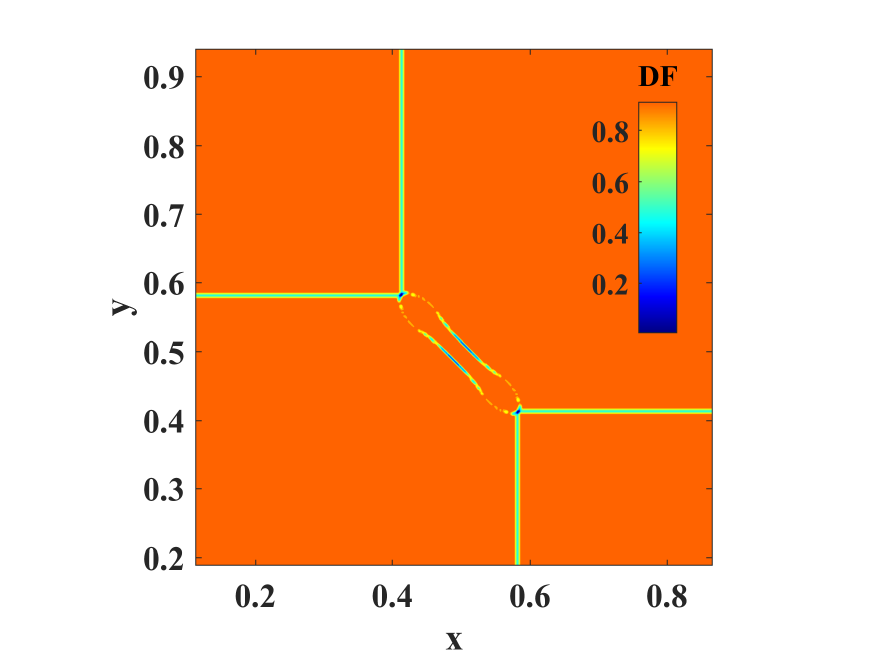}
	}\hspace{-5mm}
	\subfigure{
		\includegraphics[width=0.33\textwidth]{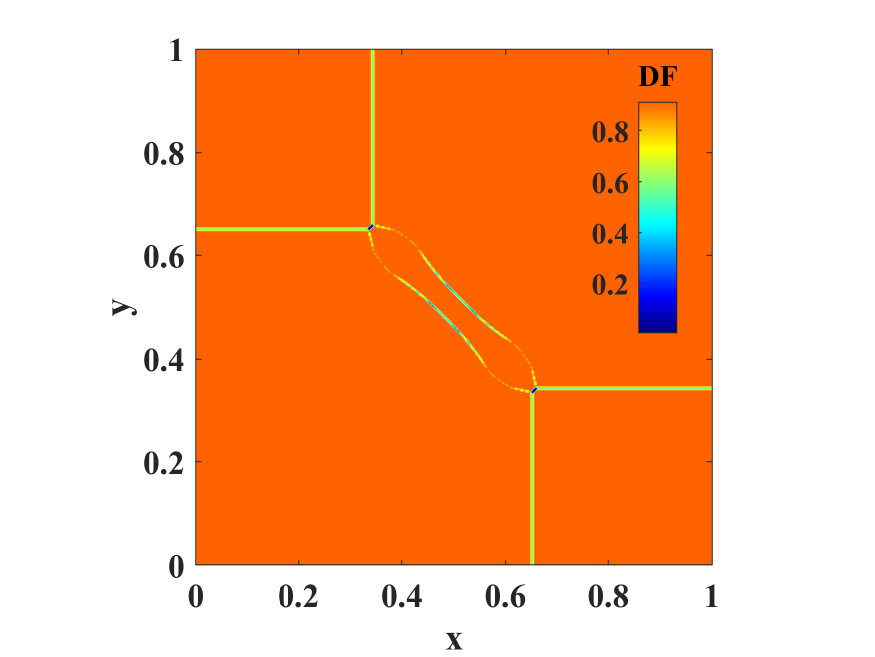}
	}\hspace{-5mm}
	\subfigure{
		\includegraphics[width=0.33\textwidth]{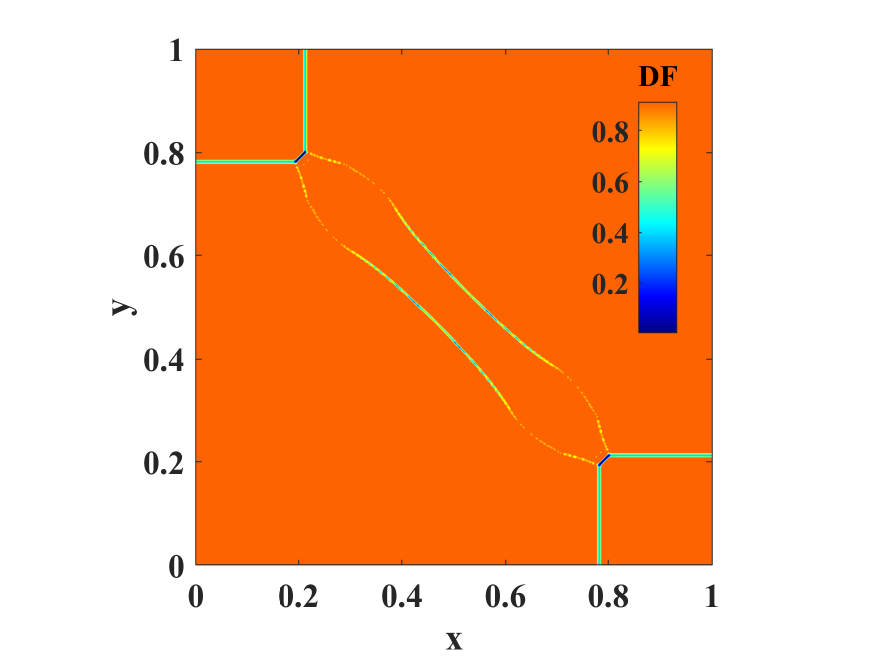}
	}
	\caption{DF distribution at step 100, 200, 400(from left to right) using S2O4 GKS solver.}
	\label{fig:27}
\end{figure}

 The solutions are shown in Fig~\ref{fig:24}-\ref{fig:27} and Table~\ref{tab:5} at $t=0.15$. The planar rarefaction waves are strong and shocks are present in the interior domain (seen in Fig~\ref{fig:24}), the hybrid reconstruction still shows good robustness in computing the extreme cases, while with limited improvement in robustness compared to the Hurricane-like problem (seen in Fig~\ref{fig:25}-\ref{fig:26}). Further, the DF approaching 0 is mainly distributed in the region where the rarefaction waves are in contact with the shocks (seen in Fig~\ref{fig:27}).

\begin{table}[htbp]
	\caption{Interaction of planar rarefaction waves: Maximum Mach number using different reconstruction methods}
	\label{tab:5}       
	\begin{tabular*}{\textwidth}{@{\extracolsep{\fill}}cccc}
		\hline\noalign{\smallskip}
		GKS solver & Maximum Mach number & LF solver & Maximum Mach number\\
		\noalign{\smallskip}\hline\noalign{\smallskip}
		WENO-AO & 1.8 & WENO-AO & 18.1\\
		Hybrid& 2.1 & Hybrid & 24.8\\
		\noalign{\smallskip}\hline
	\end{tabular*}
\end{table}


\section{Conclusion}
In this paper, we have developed a robustness-enhanced technique for the high-order finite volume scheme under two-dimensional structural mesh, and have been successfully applied it to both the high-order GKS solver and the L-F solver. This approach incorporates a newly designed hybrid spatial reconstruction based on a combination of DF and WENO methods, targeting on enhancing the robustness of the classical high-order reconstruction while maintaining high resolution. One of the challenges faced by the WENO method is its stencil selection strategy when strong discontinuities exist between the target cell and its neighboring cells. To address this issue, we have introduced a straightforward stencil selection strategy through the DF distribution to enhance the algorithm's robustness. The accuracy, high resolution, and robustness of the algorithm have been verified through numerical examples spanning from smooth flows to scenarios involving strong discontinuities. 
The current hybrid method excels at capturing complex flow structures generated by high-frequency waves, shocks, shear layers, etc., and has only exhibits slightly reduced resolution compared to the WENO method. 
In the robustness cases, DF approximate to 0 are concentrated in the region where strong shock and rarefaction waves exist and the high-order scheme will naturally reduce to the first-order one, which improves the robustness of the algorithm. 
The new idea operates primarily on spatial reconstruction, and can therefore be easily generalized to methods with arbitrary higher-order accuracy on structured or unstructured meshes. In order to further explore the potential of the new technique, future investigations will be conducted based on compact schemes like compact GKS and DG, etc.

%
%

\begin{acknowledgements}
The current research is supported by National Numerical Windtunnel project  and  National Natural Science Foundation of China (12172316, 12302378).
\end{acknowledgements}

%

\section*{Conflict of interest}
The authors declare that they have no conflict of interest.


%
%
\bibliographystyle{spmpsci}
\bibliography{reference}

\end{document}